\documentclass[12 pt]{article}

\usepackage{helvet} 

\usepackage{graphicx}  
\usepackage{amssymb}
\usepackage{amsmath}
\usepackage{subcaption}

\usepackage{makecell}
\usepackage{multirow}

\usepackage{rotating}
\usepackage{tikz}
\usetikzlibrary{arrows,
                positioning,
                shapes.geometric}

\usepackage{adjustbox}
                
\usepackage{hyperref}
\hypersetup{
    colorlinks=true,
    linkcolor=blue,
    citecolor=blue,
    filecolor= magenta,      
    urlcolor= magenta 
}

\newenvironment{keywords}{\textsf{Keywords:}\hspace{\stretch{1}}}{\hspace{\stretch{1}}\rule{1ex}{1ex}}

\newcommand{\rr}{\mathbf{r}}

\usepackage[margin=1.0in]{geometry}

\title{Efficient Generation of Membrane and Solvent Tetrahedral Meshes for Ion Channel \\ Finite Element Calculation}
\author{Zhen Chao\thanks{Department of Mathematics, University of Michigan, Ann Arbor, MI 48109, USA} \; 
\and 
Sheng Gui$^{\dag \ddagger}$ \and Benzhuo Lu\thanks{LSEC, NCMIS, Academy of Mathematics and Systems Science, Chinese Academy of Sciences, Beijing 100190, China 
} \thanks{School of Mathematical Sciences,  University of Chinese Academy of Sciences, Beijing 100049, China}\; 
\and 
Dexuan Xie\thanks{Corresponding author. E-mail address: dxie@uwm.edu. Department of Mathematical Sciences,
University of Wisconsin-Milwaukee, Milwaukee, WI 53201, USA}}

\date{January 14, 2022}	

\begin{document}


\maketitle

\begin{abstract}
A finite element solution of an ion channel dielectric continuum model such as Poisson-Boltzmann equation (PBE) and a system of Poisson-Nernst-Planck equations (PNP) requires tetrahedral meshes for an  ion channel protein region, a  membrane region, and an ionic solvent region as well as an interface fitted irregular tetrahedral mesh of a  simulation box domain. However, generating these meshes is very difficult and highly technical due to the related three regions having very complex geometrical shapes. Currently, an ion channel mesh generation software package developed in Lu's research group is one available in the public domain. To significantly improve its mesh quality and computer performance, in this paper, new numerical schemes for generating membrane and solvent meshes are presented and implemented in Python, resulting in a new  ion channel mesh generation software package.  Numerical results are then reported to demonstrate the efficiency of the new numerical schemes and the  quality of meshes generated by  the new package for ion channel proteins with ion channel pores having different geometric complexities. 
\end{abstract}







\section{Introduction}

The Poisson-Boltzmann equation (PBE) \cite{lu2008recent, holst1994poisson, chen2007finite, xie2014new} and a system of  Poisson-Nernst-Planck (PNP) equations  \cite{horng2012pnp, xie2020finitezhen, chao2021improved}  are two commonly-used dielectric continuum models for simulating an ion channel protein embedded in an membrane and immersed  in an ionic solvent. While PBE is mainly used to calculate electrostatic solvation and binding free energies, PNP is an important tool for computing membrane potentials, ionic transport fluxes, conductances,  and electric currents, etc. Both PBE and PNP have been solved approximately by typical numerical techniques such as finite difference, finite element, and  boundary element methods. Among these  techniques, finite element techniques can be more suitable to deal with the numerical difficulties caused by the complicated interface and boundary value conditions of PBE and PNP. Due to using unstructured tetrahedral meshes, they allow us to well retain  the geometry shapes of protein, membrane, and solvent regions  such that we can obtain a PBE/PNP numerical solution in a high degree of accuracy. 

However, generating an unstructured tetrahedral mesh workable for PBE/PNP finite element calculation can be very difficult and highly technical because a membrane region can cause a solvent region to have a very complicated geometrical shape, not mention that how to generate a membrane mesh remains a challenging research topic. In fact, because of the lack of membrane molecular structural data, generating a membrane mesh can become very difficult. The generation of a solvent mesh can be another major difficulty since a  solvent region can have very complex interfaces with protein and membrane regions. To avoid these difficulties,  a novel two-region approach is usually adopted to the development of an ion channel mesh generation scheme. That is, an ion channel simulation box is first divided into an ion channel protein region surrounded by an expanded solvent region without involving any membrane, where a mesh of this  expanded solvent region is supposed to have been properly constructed such that it contains both membrane and solvent meshes that are to be constructed for an ion channel PBE/PNP finite element calculation; thus, the next work to be done is to develop a numerical scheme for extracting these two membrane and solvent meshes from the  expanded solvent mesh. Based on this novel two-region approach,  an ion channel mesh software package was developed in Lu's research group with more than five-years efforts (2011 to 2017)  \cite{chen2012triangulated,tu2014software,liu2015membrane, liu2018efficient, chen2011tmsmesh, liu2017quality}. After five years,  this package remains the unique one available in the public domain and applicable for PBE/PNP ion channel finite element calculation. For clarity, we will refer it as  ICMPv1 (i.e., Ion Channel Mesh Package version 1). 


Clearly, the quality of membrane and solvent meshes extracting from an expanded solvent region strongly depends on the construction of a mesh extraction numerical scheme.  In 2014 \cite{tu2014software},  cylinders (or spheres)  were suggested to use in the separation of the membrane and pore regions but a separation process was mainly done manually. The first numerical extraction scheme  was reported in 2015 \cite{liu2015membrane}, which significantly improved the usage and performance of ICMPv1. In this scheme, a walk-detect method was adapted to detect the inner surface of an ion channel pore numerically, making it possible for us to generate the solvent, membrane, and protein region meshes and an interface fitted mesh of a simulation box domain without involving any manual effort. However, the mesh quality and the performance of ICMPv1 rely on the selection of walk step size, the number of searching layers, and the other  mesh generation parameters. A proper selection of the values of these parameters is turned out to be difficult and very time-consuming for an ion channel protein having an ion channel pore with a complicated geometrical shape.




Recently, ICMPv1 was adapted  to the implementation of the new PBE/PNP ion channel  finite element solvers developed in Xie's research group \cite{xie2020finitezhen, chao2021improved,Xie4PNPicPeriodic2020,Xie4nuPBEic2022}.  During these applications,  ICMPv1 was found to occasionally produce a membrane  mesh that contains the tetrahedra belonging to a solvent mesh. It is possible to remove these false tetrahedra manually. For example, via a visualization tool (e.g.,  ParaView \cite{ahrens2005paraview}), we may  identify these false tetrahedra and then remove them by reconstructing a new mesh. We may also adjust the related parameters repeatedly until none of  false tetrahedra occur in a membrane mesh, but doing so may be very time-consuming and may  twist a protein, membrane, or solvent mesh due to using improper parameter values, causing the numerical accuracy of a PBE/PNP finite element solution to be reduced significantly. These cases motivated Xie's research group to develop more effective and more efficient numerical schemes than those used in ICMPv1. Eventually,  the second version of ICMPv1, denoted by ICMPv2, has been developed by Xie's research group through a close collaboration with Lu's research group. The purpose of this paper is to present the new schemes implemented in ICMPv2 and report numerical results to demonstrate that ICMPv2 can generate membrane and solvent meshes much more efficiently and  in much higher quality than ICMPv1, especially for an ion channel protein having an ion channel pore with a complicated geometric shape even for an ion channel protein with multiple ion channel pores --- a case in which ICMPv1 would not work.


The rest of the paper is organized as follows. Section 2 introduces a framework for ion channel  mesh generations. Section 3 presents a new scheme for generating a triangle surface mesh of the box domain boundary.  Section 4 presents a new scheme for selecting mesh points to be set as the vertices of an interface triangular mesh between the membrane and solvent meshes. Section 5 presents a new numerical scheme for extracting the membrane and solvent tetrahedral meshes from an expanded solvent tetrahedral mesh. Section 6 reports numerical results to demonstrate that ICMPv2 can generate membrane and solvent meshes in higher quality and better computer performance  than ICMPv1. Conclusions are made in Section~7.

\section{A framework for ion channel mesh generation}
\begin{figure}[t]
 \centering
         \begin{minipage}[b]{0.45\textwidth}
                \centering
                \includegraphics[width=\textwidth]{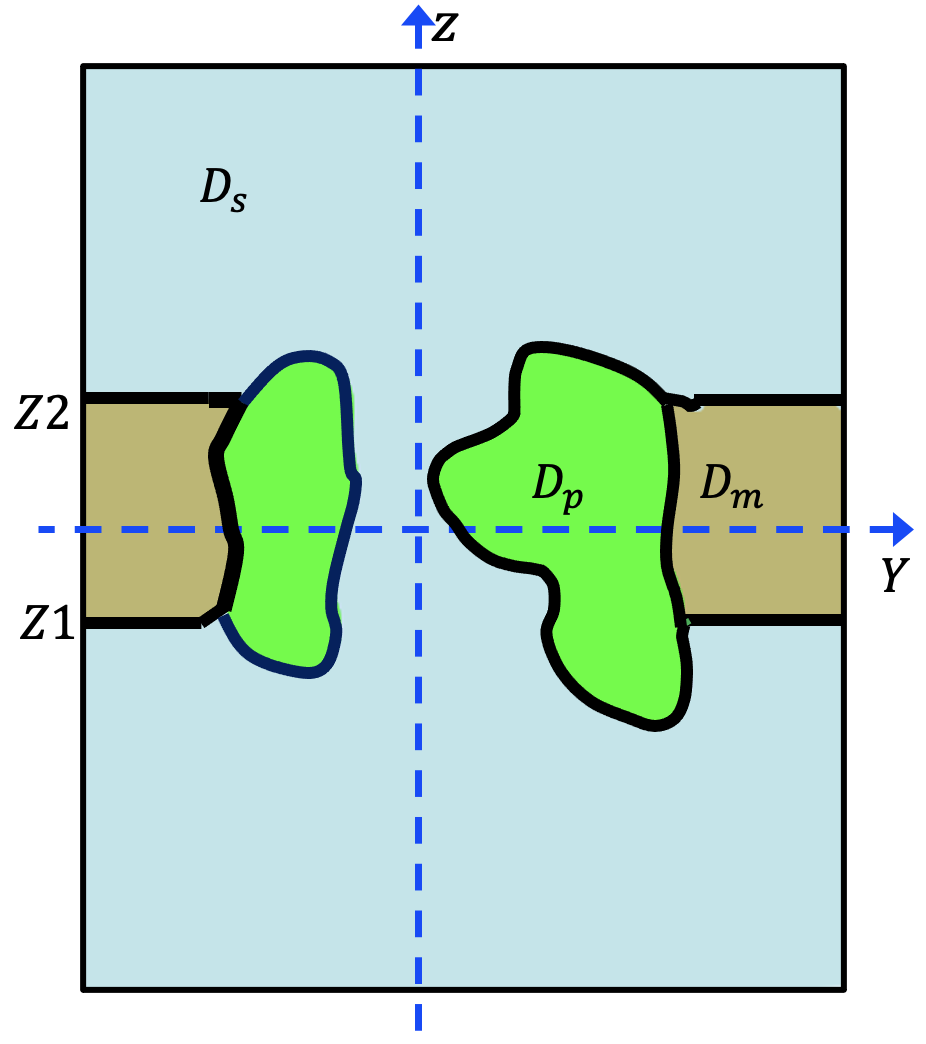}
                \caption{An illustration of box domain partition~\eqref{DomainPartition}. }
                \label{fig:3Ddomain}
        \end{minipage}
        \qquad
        \begin{minipage}[b]{0.45\textwidth}
                \centering
                \includegraphics[width=\textwidth]{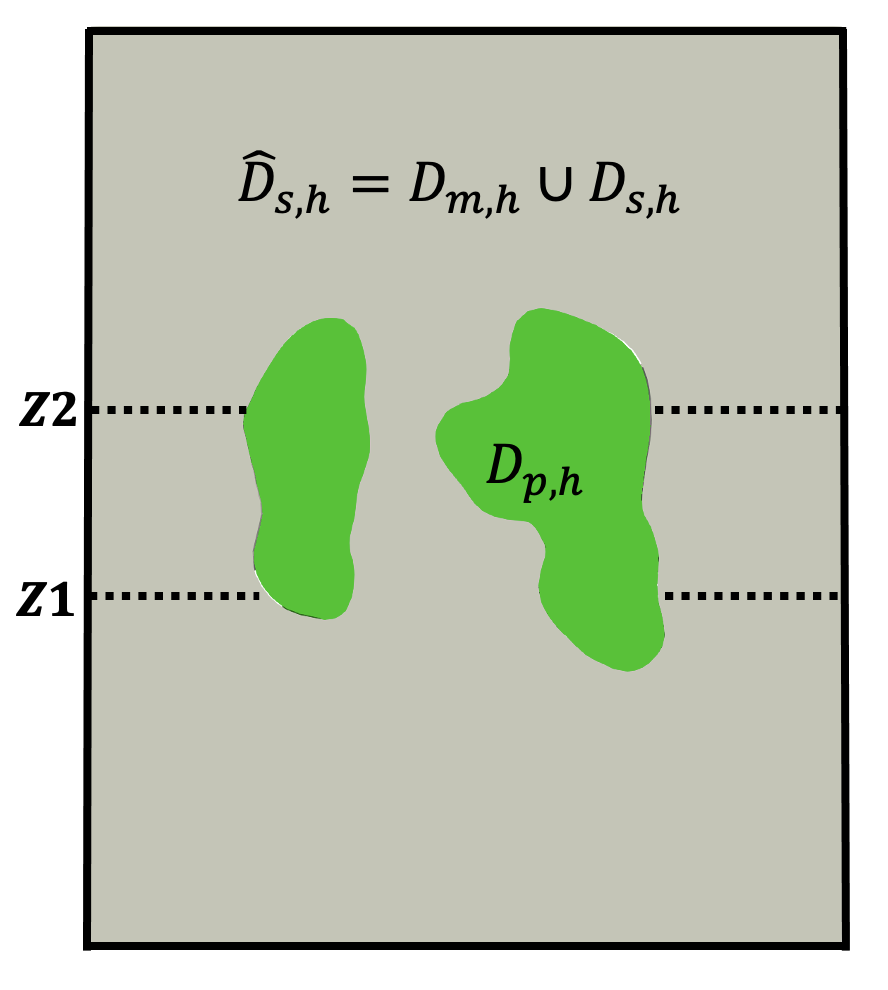}
                \caption{An illustration of the mesh partitions of  \eqref{meshDomain2}.   }        
                 \label{twoRegionPartition}
        \end{minipage}  
\end{figure}

As required to implement a continuum dielectric  model such as PBE or PNP for an ion channel system  consisting of an ion channel protein, a membrane, and an ionic solvent, a sufficiently large simulation box, $\Omega$, is selected to satisfy the domain partition
\begin{equation}\label{DomainPartition}
\Omega = {D}_p \cup {D}_m \cup {D}_s,
\end{equation}
where $ {D}_p$, $ {D}_m$, and $ {D}_s$  denote a protein region, a membrane region, and a solvent region, respectively.  Specifically,  we construct $\Omega$ by
\begin{equation}\label{meshDomain}
\Omega= \left\{(x, y, z) | L_{x_1} < x  <  L_{x_2}, L_{y_1}  <  y  <  L_{y_2}, L_{z_1}  <  z <  L_{z_2}\right\},
\end{equation}
where $L_{x_1}$, $L_{x_2}$, $L_{y_1}$, $L_{y_2}$, $L_{z_1}$, and $L_{z_2}$ are real numbers. We then assume that the center of an ion channel pore is at the origin of the rectangular coordinator system and the location of $D_m$ is determined by two parameters $Z1$ and $Z2$ between $L_{z_1}$ and $L_{z_2}$. An illustration of our box domain construction is given in Figure \ref{fig:3Ddomain}.

However, in practice, the subdomains $D_p$, $D_m$, and $D_s$ are approximated by  their tetrahedral meshes ${D}_{p,h}$, $ {D}_{m,h}$, and $ {D}_{s,h}$ since it is too difficult to obtain them directly due to their complex geometric shapes. As soon as ${D}_{p,h}$, $ {D}_{m,h}$, and $ {D}_{s,h}$ are obtained, an interface fitted tetrahedral mesh, $\Omega_h$, of $\Omega$ can be constructed according to  the following mesh domain  partition
\begin{equation}\label{meshDomain1}
\Omega_h = {D}_{p,h} \cup {D}_{m,h} \cup {D}_{s,h}.
\end{equation} 

One key step to obtain ${D}_{p,h}$, $ {D}_{m,h}$, and $ {D}_{s,h}$   is to construct their  triangular surface meshes  $\partial{D}_{p, h}$, $\partial{D}_{m, h}$, and $\partial{D}_{s, h}$ since  from these  three triangular surface meshes the tetrahedral volume meshes ${D}_{p,h}$, $ {D}_{m,h}$, and $ {D}_{s,h}$ can be generated routinely by using a volume mesh generation software package such as TetGen \cite{si2005meshing, si2015tetgen}. 

It has been known that  $\partial{D}_{p, h}$ can be generated from one of the current software packages TMSmesh \cite{chen2012triangulated},  NanoShaper \cite{decherchi2013general},  GAMer \cite{yu2008high}, MSMS \cite{sanner1996reduced},  and MolSurf \cite{sjoberg1997molsurf} when a molecular structure of an ion channel is known. However, how to  generate $\partial{D}_{m,h}$ remains a difficult research topic. In fact, a membrane consists of  a double layer of phospholipid, cholesterol, and glycolipid molecules, making it very difficult to derive a boundary mesh of  $\partial{D}_{m,h}$. 

To avoid these difficulties,  we follow the strategy used in \cite{tu2014software,liu2015membrane} to construct  an expanded solvent tetrahedral mesh, $ \hat{D}_{s,h}$,  satisfying 
\begin{equation}\label{meshDomain2}
 \Omega_h = D_{p,h} \cup \hat{D}_{s,h}, \quad \hat{D}_{s,h} = D_{m,h} \cup D_{s,h}.
\end{equation}
An illustration of the above partition is given in Figure~\ref{twoRegionPartition}, where the dotted lines represents a set of mesh points,  denoted by ${\cal S}$, to be applied to the construction of $ \hat{D}_{s,h}$.  A numerical scheme for a selection of ${\cal S}$ will be presented in Section~4  to ensure that  $\hat{D}_{s,h} $ contains both $D_{m,h}$ and $D_{s,h}$.

 Clearly, the boundary surface mesh $\partial \hat{D}_{s,h}$ of $ \hat{D}_{s,h}$ can be constructed by
\[     \partial \hat{D}_{s,h} = \partial D_{p,h} \cup \partial \Omega_h,\]
where $ \partial D_{p,h}$ has been given and  $\partial \Omega_h$ is a triangular surface mesh of the boundary $\partial \Omega$ of the box domain $\Omega$, whose construction will be done by a numerical scheme to be presented in Section~3.  We then can use current mesh software packages (e.g., TetGen) to generate the tetrahedral meshes $D_{p,h}$ and $ \hat{D}_{s,h}$.  As soon as $\hat{D}_{s,h}$ is known,  we  need a numerical scheme to extract ${D}_{m,h}$ and ${D}_{s,h}$ from $\hat{D}_{s,h} $. Such a scheme will be presented in Section~5. 


\section{Construction of a box triangular surface mesh}
In this section, we present a numerical scheme for constructing a triangular surface mesh, $\partial \Omega_h$, to ensure that the mesh domain partition \eqref{meshDomain1} holds. In this scheme, a triangular surface mesh,  $\partial D_{p,h}$, is supposed to be given. Thus, we can find the  smallest rectangular box $[a_1, b_1]  \times[a_2, b_2]\times[a_3, b_3]$ that holds   $\partial D_{p,h}$ using the  formulas:
\begin{equation}\label{proteinSize}
\begin{aligned}
a_1&= \min_{1\leq i\leq N} x_i,\qquad b_1 =  \max_{1\leq i\leq N} x_i, \\
a_2&=  \min_{1\leq i\leq N} y_i, \qquad  b_2=  \max_{1\leq i\leq N} y_i, \\
a_3&=  \min_{1\leq i\leq N} z_i,  \qquad  b_3=  \max_{1\leq i\leq N} z_i, 
\end{aligned}
\end{equation}
where ($x_i$, $y_i$, $z_i$) denotes the position vector of the $i$-th mesh point of $\partial D_{p,h}$ and $N$ is the total number of  mesh points of $\partial D_{p,h}$. We then can construct a box domain, $\Omega$, of \eqref{meshDomain}  in terms of  three parameters, $\eta_i$ for $i=1,2,3$,  according to the following formulas:
 \begin{equation}\label{eq:boxDomainSize}
 \begin{aligned}
                  L_{x_1} = a_1 - \eta_1, \qquad L_{x_2} = b_1+ \eta_1,\\
	           L_{y_1} = a_2- \eta_2,\qquad  L_{y_z} = b_2+ \eta_2,\\
	              L_{z_1} = a_3 - \eta_3, \qquad  L_{z_2} =b_3 + \eta_3.
\end{aligned}	              
\end{equation}
The default values of $\eta_1, \eta_2$, and $\eta_3$ are set as 20 but can be adjusted by users as needed. In this way, a selection of a box domain satisfying the partition \eqref{DomainPartition} is greatly simplified. 

\begin{figure}[t]
 \centering
        \begin{minipage}[b]{0.45\textwidth}
                \centering
                \includegraphics[width=0.73\textwidth]{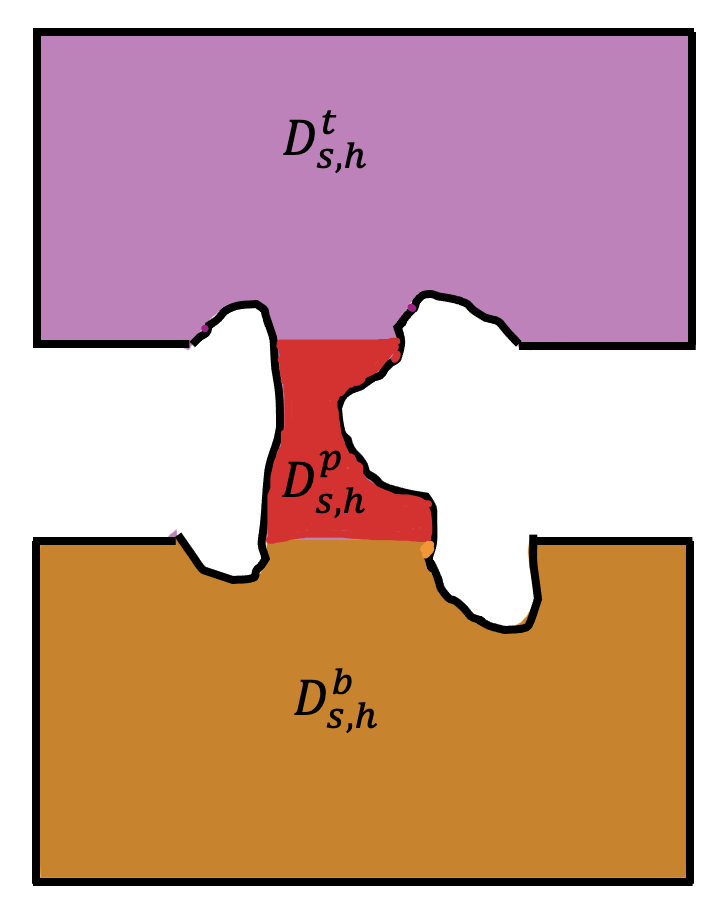}
                \caption{An illustration of the partition \eqref{solvenDomainPartition} of a solvent mesh $D_{s,h}$.}
                \label{fig:1magDsPartition}
       \end{minipage} 
\qquad
         \begin{minipage}[b]{0.45\textwidth}
                \centering
                \includegraphics[width=\textwidth]{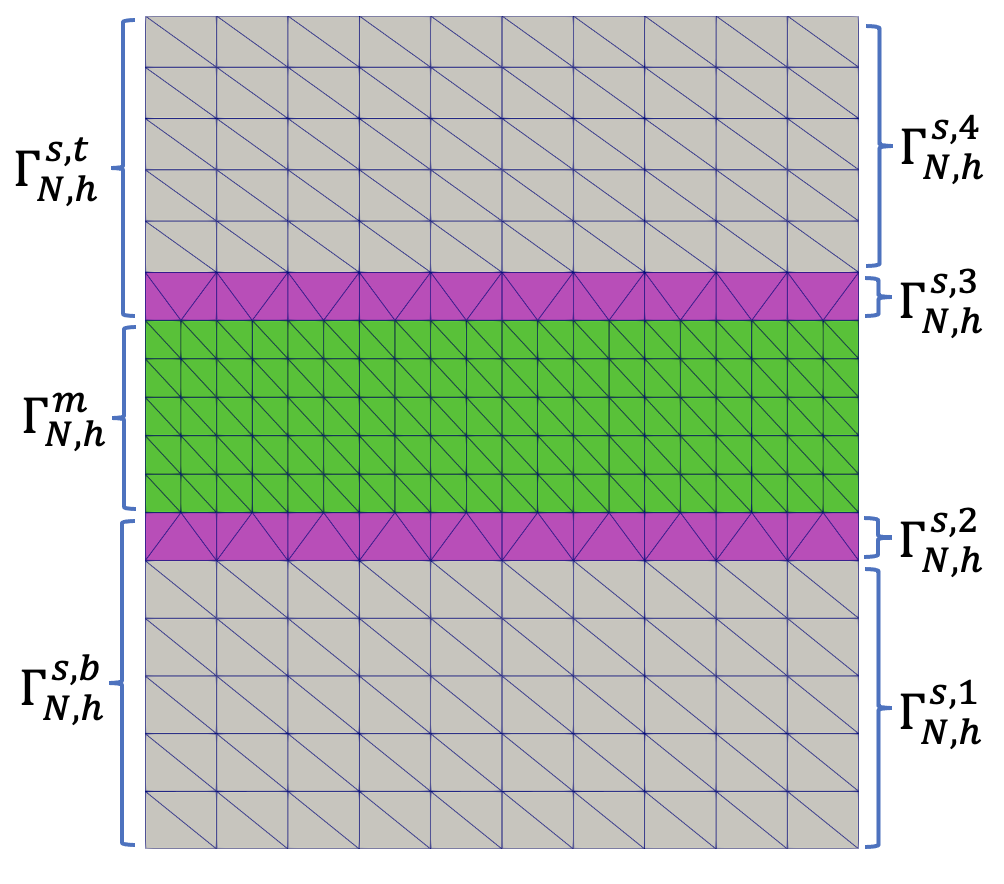}
                \caption{An example of our box surface mesh partitions \eqref{lateralPartition} and \eqref{layerPartition}. }
                \label{fig:boundaryPartition}
        \end{minipage}        
\end{figure}

We next split the six boundary surfaces of  $\partial\Omega$  by
\begin{equation}\label{eqPartition}
 \partial\Omega= \Gamma_D \cup  \Gamma_N, 
\end{equation}
and divide  a solvent mesh, $D_{s,h}$,   into three portions, $D_{s,h}^b, D_{s,h}^t,$ and $D_{s,h}^p$, by
\begin{equation}\label{solvenDomainPartition}
D_{s,h} = D_{s,h}^b \cup D_{s,h}^t \cup D_{s,h}^p,
\end{equation}
where   $\Gamma_{D}$ and $\Gamma_N$ contain the bottom and top surfaces and the four side surfaces of $\partial\Omega$,  respectively, and $D_{s,h}^b$,  $D_{s,h}^t$, and  $D_{s,h}^p$ are defined by
\[ D_{s,h}^b=\{ \rr \in D_{s,h}  \; | \;  \rr =(x,y,z) \mbox{ with } z < Z1\}, \]
\[ D_{s,h}^t=\{ \rr \in D_{s,h}  \; | \; \rr =(x,y,z) \mbox{ with } z > Z2 \},\]
and 
\[ D_{s,h}^p=\{ \rr \in D_{s,h}  \;  | \; \rr =(x,y,z) \mbox{ with } Z1 \leq z \leq Z2\}.\]
An illustration of partition \eqref{solvenDomainPartition} is given in Figure~\ref{fig:1magDsPartition}.

Since a uniform triangular mesh of $\Gamma_{D}$ can be constructed easily, we only describe the construction of a triangular surface mesh, $\Gamma_{N,h}$, of $\Gamma_{N}$.  According to the solvent mesh partition  \eqref{solvenDomainPartition}, we can split   $\Gamma_{N,h}$ into three sub-meshes, $\Gamma_{N,h}^{s,b}$, $ \Gamma_{N,h}^{m}$, and $\Gamma_{N,h}^{s,t}$, such that
\begin{equation}\label{lateralPartition}
\Gamma_{N,h} = \Gamma_{N,h}^{s,b} \cup \Gamma_{N,h}^{m}\cup \Gamma_{N,h}^{s,t}, 
\end{equation}
where $\Gamma_{N,h}^{s,b}$, $ \Gamma_{N,h}^{m}$, and $\Gamma_{N,h}^{s,t}$ are defined by
\[ \Gamma_{N,h}^{s,b} = \Gamma_{N}\cap D_{s,h}^b, \; \Gamma_{N,h}^{s,t} = \Gamma_{N}\cap D_{s,h}^t, \; \Gamma_{N,h}^{m} = \Gamma_{N}\cap D_{m,h}.\]

In ICMPv2, these three sub-meshes are constructed as uniform triangular meshes, respectively, and are allowed to have  different mesh sizes. In particular, we set the mesh size $h_m$ of $\Gamma_{N,h}^{m}$ as an input parameter of  ICMPv2. We then set $\Gamma_{N,h}^{s,b} $ and $\Gamma_{N,h}^{s,t} $ to have the same mesh size $h_s$. By default, we set $h_s=2h_m$. In this case,  $\Gamma_{N,h}^{s,b} $ and  $\Gamma_{N,h}^{s,t}$ can be split as follows:
 \begin{equation}\label{layerPartition}
 \Gamma_{N,h}^{s,b}  = \Gamma_{N,h}^{s,1} \cup \Gamma_{N,h}^{s,2},\quad \Gamma_{N,h}^{s,t}  = \Gamma_{N,h}^{s,3} \cup \Gamma_{N,h}^{s,4},
 \end{equation} 
 where $\Gamma_{N,h}^{s,2}$ and $\Gamma_{N,h}^{s,3}$ are defined by one mesh layer of $\Gamma_{N,h}^{s,b}$ and $\Gamma_{N,h}^{s,t}$, respectively. An illustration of surface mesh partitions \eqref{lateralPartition} and \eqref{layerPartition} is given in Figure~\ref{fig:boundaryPartition}, where $\Gamma_{N,h}^{s,2}$ and $\Gamma_{N,h}^{s,3}$  are  colored in pink, $\Gamma_{N,h}^{s,1}$ and $\Gamma_{N,h}^{s,4}$ in grey, and $\Gamma_{N, h}^m$  in green. The reason why we set $h_m = h_s/2$ is to select a sufficiently large set of ${\cal S}$ from the membrane top and bottom surfaces to improve the numerical quality of a membrane mesh, $D_{m,h}$, to be extracted from an expanded solvent mesh, $\hat{D}_{s,h}$.  
 
%
 
\section{A numerical scheme for selecting membrane surface mesh points}
In this section, we present a numerical scheme for selecting a membrane surface mesh point set, ${\cal S}$, as needed in the construction of an expanded solvent mesh, $\hat{D}_{s,h}$, satisfying \eqref{meshDomain2}. For clarity, we only describe a mesh point selection from the bottom membrane surface $\Gamma_{m}^{b}$ since a selection from the top membrane surface can be done similarly.


\begin{figure}[t]
 \centering
        \begin{subfigure}[t]{0.45\textwidth}
                \centering
                \includegraphics[width=\linewidth]{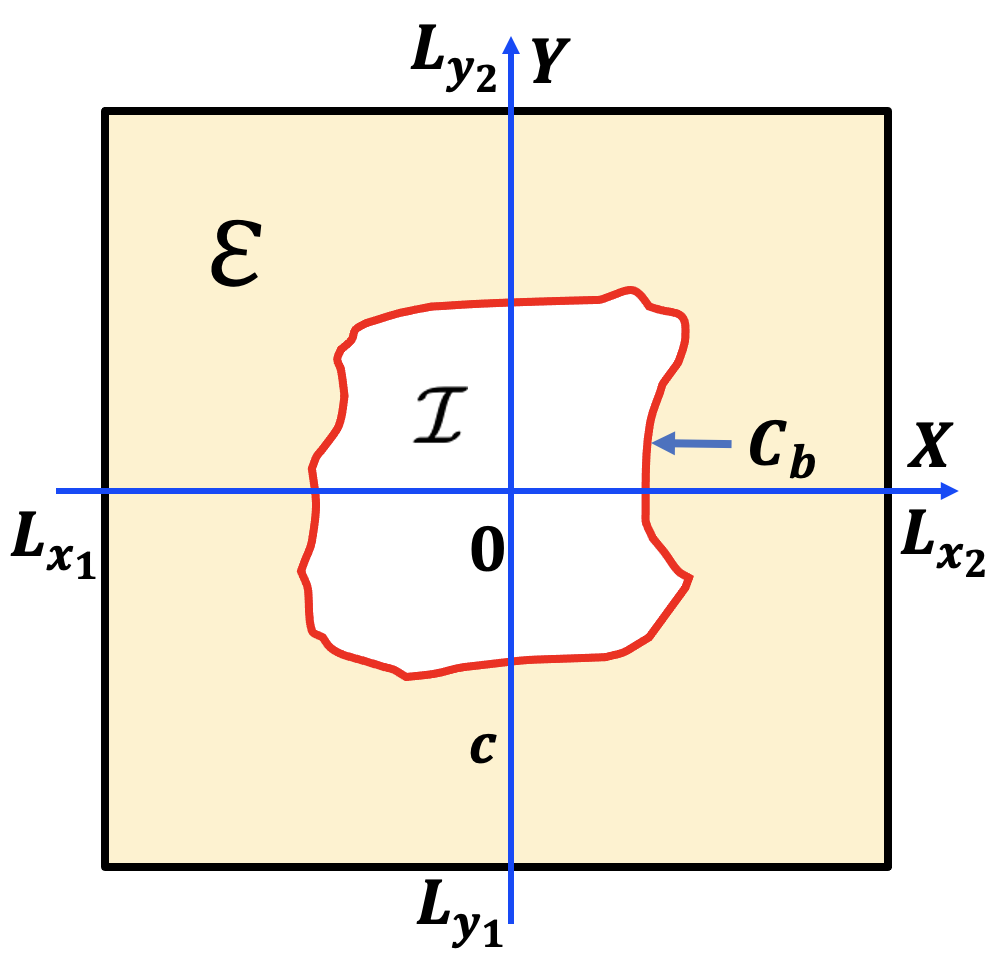}
                \caption{A partition of rectangle $[L_{x_1}, L_{x_2}]\times [L_{y_1},L_{y_2}]$ into an exterior  region, ${\cal E}$, and an interior region, ${\cal I}$. }
                \label{fig:planeAtza}
       \end{subfigure} 
\qquad
         \begin{subfigure}[t]{0.45\textwidth}
                \centering
                \includegraphics[width=\textwidth]{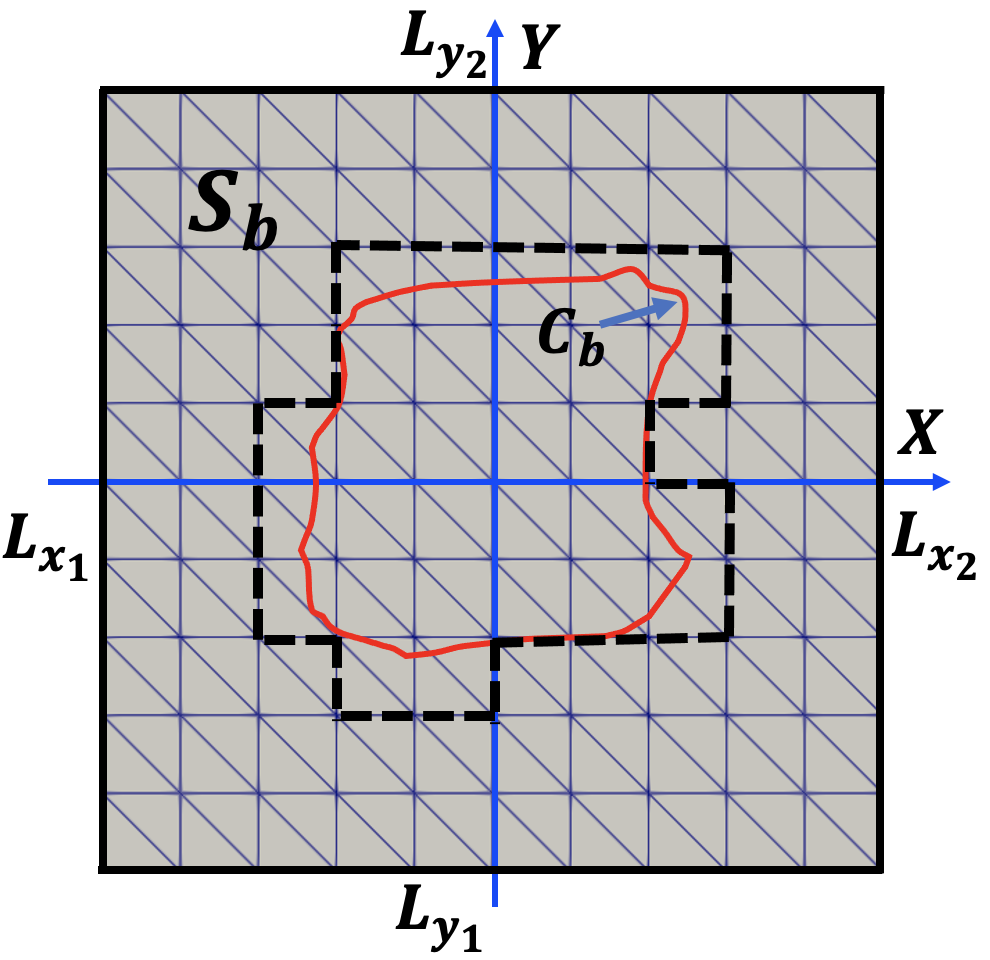}
                \caption{$S_b$ consists of the mesh points outside the dash line plus the mesh points on the dash line. }
                \label{fig:planeAtzb}
        \end{subfigure}  
         \caption{ An illustration of  a partition of rectangle $[L_{x_1}, L_{x_2}]\times [L_{y_1},L_{y_2}]$ by  a cross section curve of  $\partial D_{p,h}$ with the plan $z=Z1$ and a mesh  point set  $S_b$ defined in \eqref{setPartition}.}
          \label{fig:planeAtz}
\end{figure}

With a given mesh size, $h_m$, of a a membrane boundary mesh $\partial D_{m,h}$,  we construct a uniform rectangular mesh of the rectangle $[L_{x_1}, L_{x_2}]\times [L_{y_1},L_{y_2}]$ and define a set,  ${\cal T}$,  of its mesh points by
\[ {\cal T} = \{ (x_i, y_j) \; | \;  x_i = L_{x_1} + i h_m,   y_j = L_{y_1} + j h_m \mbox{ for } i=0,1,\ldots,m, j=0,1,\ldots,n \},\]
where $m=(L_{x_2} - L_{x_1})/h_m$ and $n=(L_{y_2} - L_{y_1})/h_m$.
 We also obtain a curve, ${\cal C}_b$, of the cross section of $\partial D_{p,h}$ with the plan $z=Z1$. Clearly, this curve splits the rectangle $[L_{x_1}, L_{x_2}]\times [L_{y_1},L_{y_2}]$ into an exterior  region, denoted by ${\cal E}$,  and  an interior   region, denoted by ${\cal I}$, as illustrated in Figure~\ref{fig:planeAtza}. We then can obtain a set, $ S_b$, of mesh points from the bottom membrane surface  by
\begin{equation}\label{setPartition}
S_b= \{(x_i, y_i, Z1) \; | \: (x_i,y_i) \in  {\cal T} \cap {\cal E} \}
\end{equation}
 as illustrated in Figure~\ref{fig:planeAtzb}.

Similarly, we can obtain a set, $S_t$,  of mesh points from the top membrane surface at $z=Z2$.  A union of $S_b$ and $S_t$ gives the set $S$ to be used in the construction of $\hat{D}_{s,h}$.

\section{Extraction of membrane and solvent meshes}
In this section, we present a numerical scheme for extracting membrane and solvent meshes, $D_{m, h}$ and ${D}_{s, h}$, from an expanded solvent mesh,  $\hat{D}_{s, h}$. In this scheme, we assume that the tetrahedra  of $D_{p, h}$ and $\hat{D}_{s, h}$ have label numbers 1 and 2, respectively, while  the tetrahedra of ${D}_{p, h}$, $D_{s, h}$, and ${D}_{m, h}$ have label numbers 1, 2, and 3, respectively. 

For clarity, we describe our new numerical scheme in six steps as follows: 

\begin{description}
\item[Step 1] Construct a rectangle, $[a, b]\times [c,d]$, that contains a portion, ${\cal P}$, of  triangular surface mesh $\partial D_{p,h}$ intercepted by the two planes $Z=Z1$ and $Z=Z2$ by 
\[
\begin{aligned}
  a &= \min\{ x_i,\; | \; (x_i, y_i, z_i)  \mbox{ on } {\cal P} \} - \tau, \; b= \max\{ x_i,\; | \;  (x_i, y_i, z_i)  \mbox{ on }  {\cal P} \} + \tau, \\
  c &=  \min\{ y_i,\; | \;  (x_i, y_i, z_i)  \mbox{ on }  {\cal P} \} - \tau,  \;  d= \max\{ y_i,\; | \;  (x_i, y_i, z_i)  \mbox{ on } {\cal P} \} + \tau, 
\end{aligned}  
  \]
 where $(x_i, y_i, z_i)$ denotes the $i$th vertex of  ${\cal P}$ and $\tau$ is a positive parameter to ensure that  the rectangle $[a, b]\times [c,d]$ is  not to touch any part of ${\cal P}$. By default, we set $\tau = h_m$. 

\begin{figure}[t]
 \centering
 \begin{minipage}[b]{0.45\textwidth}
                \centering
                \includegraphics[width=0.85\textwidth]{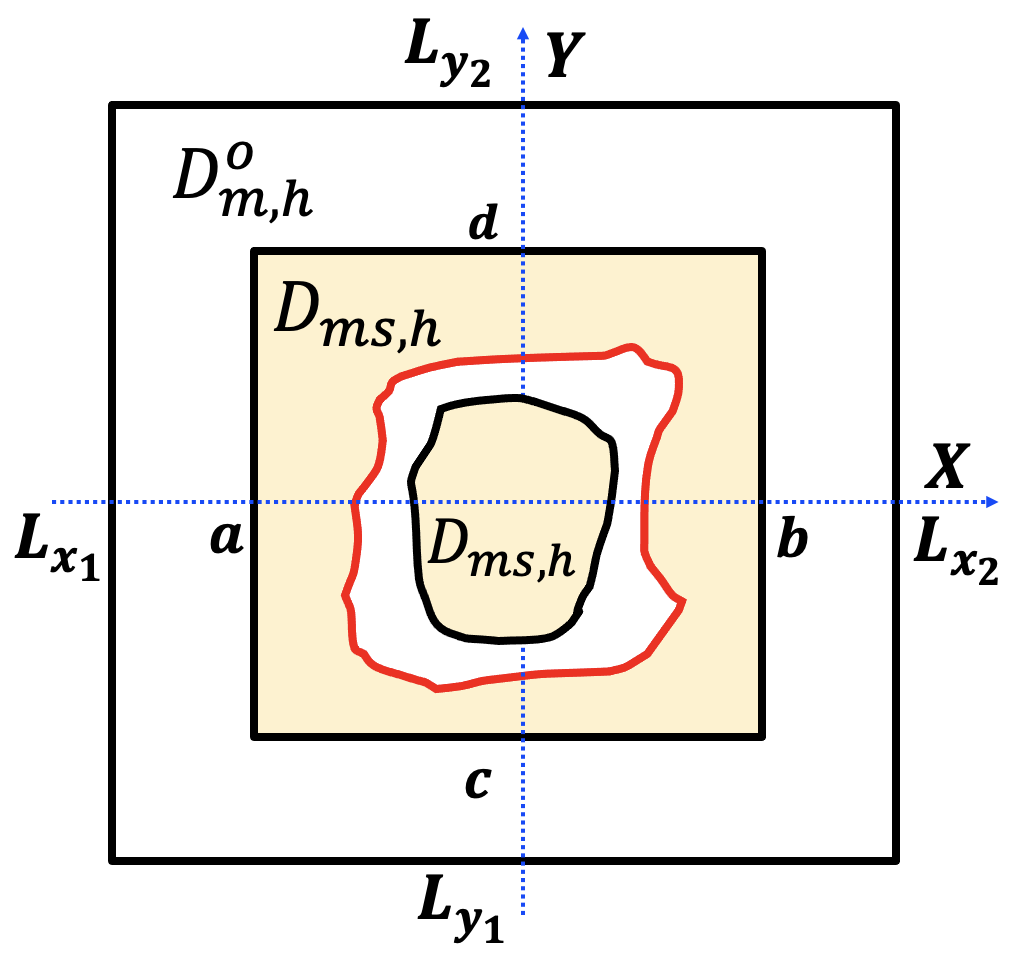}
                \caption{A partition \eqref{DmsPartition}.}
\label{fig:planeZPartition1}
 \end{minipage} 
\qquad
         \begin{minipage}[b]{0.45\textwidth}
                \centering
                \includegraphics[width=0.85\textwidth]{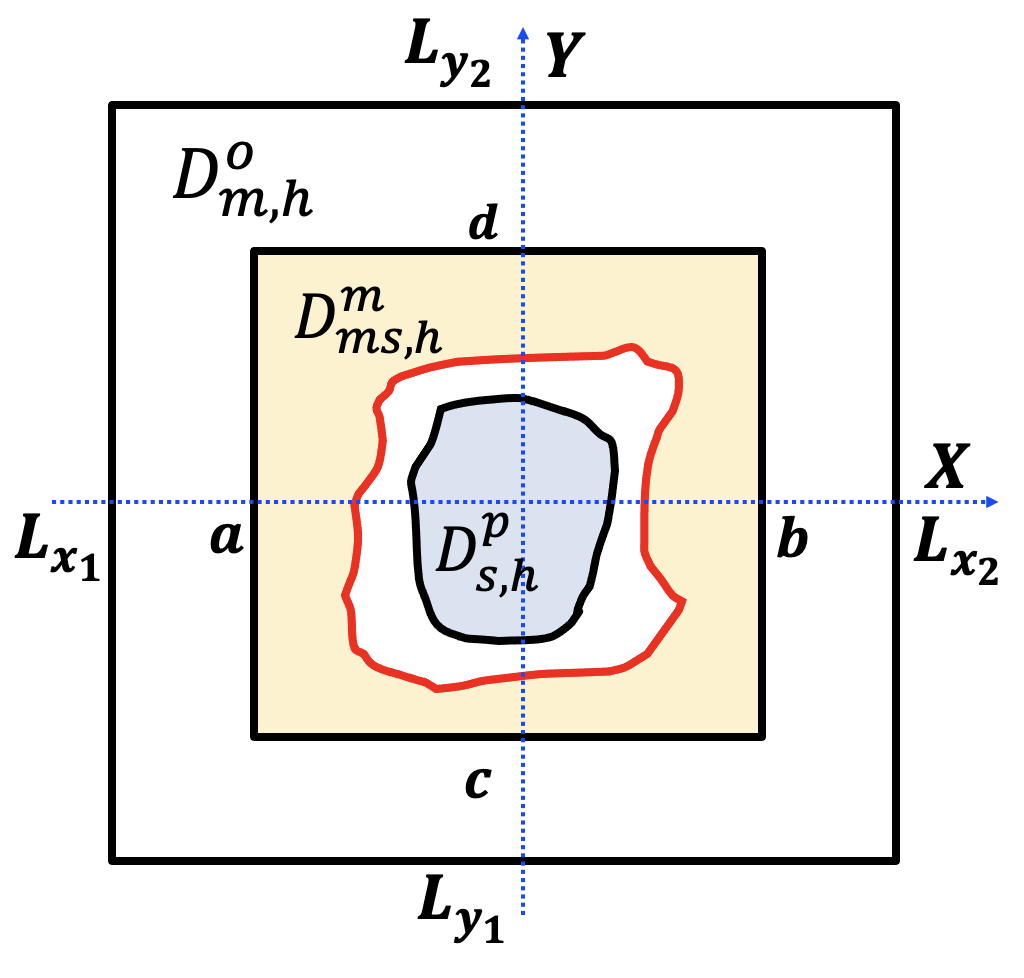}
                \caption{A partition \eqref{submesh_Ds2}.}
                \label{fig:planeZPartition2}   
        \end{minipage}    
  \begin{minipage}[b]{\textwidth}
                \centering
     \begin{subfigure}[b]{0.32\textwidth}
                \centering
               \includegraphics[width=0.8\textwidth]{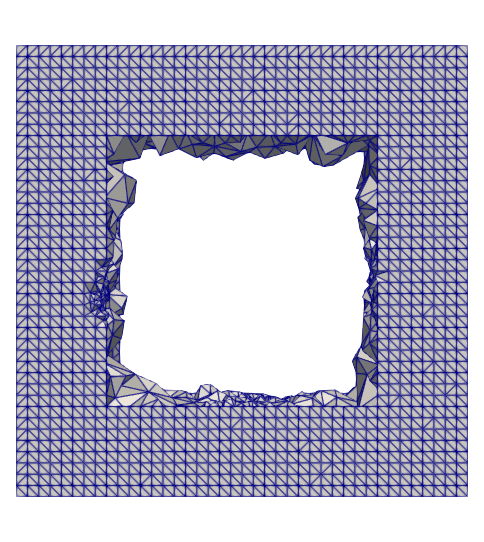}
                \caption{${D}_{m,h}^o$}
               \label{fig:3emnMembraneSubmesh}
        \end{subfigure}  
         \begin{subfigure}[b]{0.32\textwidth}
                \centering
                \includegraphics[width=1.2\textwidth]{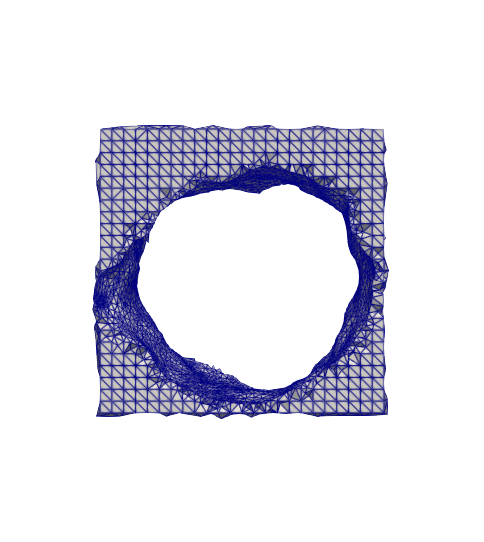}
                \vspace{-1.4cm}
                \caption{${D}_{ms,h}^m$}
               \label{fig:3emnMembraneBoundarymesh}
        \end{subfigure}  
        \;\;
         \begin{subfigure}[b]{0.32\textwidth}
                \centering
                \includegraphics[width=0.94\textwidth]{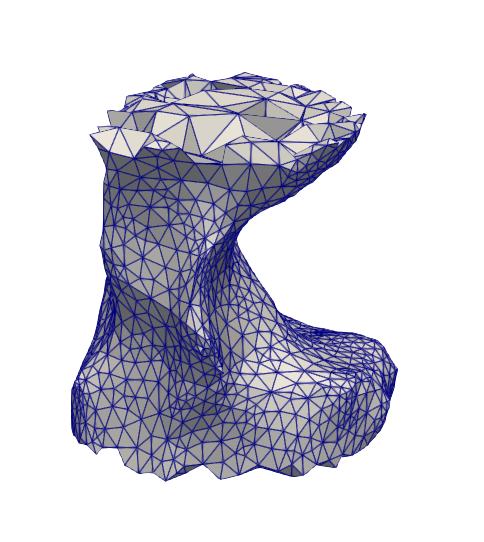}
                \vspace{-5mm}
                \caption{$ D_{s,h}^p $}
              \label{fig:3emnPoreBoundarymesh}
        \end{subfigure}   
    \caption{Submeshes ${D}_{m,h}^o$,  $D_{ms,h}^m$ and $ D_{s,h}^p $ generated by ICMPv2 for VDAC (PDB ID: 3EMN).}
                \label{fig:3emnPoreAndMembraneBoundarymesh}   
        \end{minipage}        
\end{figure} 

\item[Step 2] Construct  three submeshes, ${\cal B}_{h}$, ${D}_{ms,h}$, of $\hat{D}_{s,h}$, and  $ {D}_{m,h}^o$   by
 \begin{equation}
 \label{submesh_Ds}
 \begin{aligned}
    &  {\cal B}_{h} =  \hat{D}_{s,h} \cap  [L_{x_1}, L_{x_2}]\times [L_{y_1},L_{y_2}] \times [Z1, Z2],& \\
      &   {D}_{ms,h} =  \hat{D}_{s,h} \cap  [a, b]\times [c,d] \times [Z1, Z2],&
\end{aligned}
 \end{equation}
 and
 \begin{equation}\label{DmsPartition}
{D}_{m,h}^o =  {\cal B}_{h} -  {D}_{ms,h}.
 \end{equation}
 In other words, we have divided $ {\cal B}_{h}$ into two submeshes satisfying 
 \[       {\cal B}_{h} =  {D}_{ms,h} \cup {D}_{m,h}^o \]
 as illustrated in Figure~\ref{fig:planeZPartition1}.  From this figure we can see that ${D}_{ms,h}$ consists of two non-overlapped parts --- one part is nothing but  the solvent portion $D_{s,h}^p$ of  the solvent mesh partition \eqref{solvenDomainPartition} and the other part, denoted by $D_{ms,h}^m$, belongs to the membrane region $D_m$. Thus, ${D}_{ms,h}$  can be expressed as 
  \begin{equation}
 \label{submesh_Ds2}
      {D}_{ms,h} =  D_{ms,h}^m  \cup D_{s,h}^p.
 \end{equation}
 As examples, the submeshes ${D}_{m,h}^o$, $D_{ms,h}^m$, and $ D_{s,h}^p $ generated by ICMPv2 for an ion channel protein (VDAC) are displayed in Figure~\ref{fig:3emnPoreAndMembraneBoundarymesh}.

 \item[Step 3] Separate the tetrahedra of  ${D}_{ms,h}$ as two non-overlapped sets, one set leads to $D_{ms,h}^m$ and the other set to $D_{s,h}^p$, such that the partition \eqref{submesh_Ds2} holds.
In ICMPv2,  this separation is done by a numerical scheme implemented in the python function {\tt split()} from the  Python library {\em Trimesh}\footnote{\em https://github.com/mikedh/trimesh}. To do so, we need to obtain a boundary triangular mesh $\partial{D}_{ms,h}$ of ${D}_{ms,h}$ since it is a required input mesh of  this python function. We obtain  $\partial{D}_{ms,h}$ through finding the boundary surface meshes of $D_{ms,h}^m$ and $D_{s,h}^p$, respectively.

\item[Step 4] Identify the tetrahedra of $D_{s,h}^p$  by doing ray tests via the  ray-triangle intersection method (see \cite{moller1997fast} for example). To do so, we first calculate the centroids of all the tetrahedra of ${D}_{ms,h}$ and then use ray tests to check if they are inside the volume region enclosed by a boundary triangular surface mesh  of $D_{s,h}^p$ or not. In the true case, we store the tetrahedron indices to an index set, ${\cal S}_s$; else, the tetrahedron indices are stored  to another index set, ${\cal S}_m$. We then change the label numbers of the tetrahedra with indices in ${\cal S}_m$ from 2 to 3 to obtain the first part of $D_{m,h}$, which is denoted by $D_{m,h}^1$.
In ICMPv2, a ray test is done by calling the python function {\tt contains\_points()} from a class object,
  {\tt   \indent trimesh.ray.ray\_pyembree.RayMeshIntersector()}
  of the  Python library {\em Trimesh}.

\item[Step 5]  Identify the tetrahedra of $ {D}_{m,h}^o$ by using partition \eqref{DmsPartition}  and change their label numbers from 2 to 3 to obtain the second part of $D_{m,h}$, which is denoted by $D_{m,h}^2$. 

\item[Step 6] Obtain the membrane and solvent meshes ${D}_{m,h}$ and ${D}_{s,h}$ by 
\[{D}_{m,h} = {D}_{m,h}^1 \cup D_{m,h}^2, \quad  {D}_{s,h} = \hat{D}_{s,h} - {D}_{m,h}. \]
\end{description}

\begin{figure}[t]
        \centering
        \begin{minipage}[b]{\textwidth}
                 \begin{subfigure}[b]{0.24\textwidth}
                \centering
                \includegraphics[width=\textwidth]{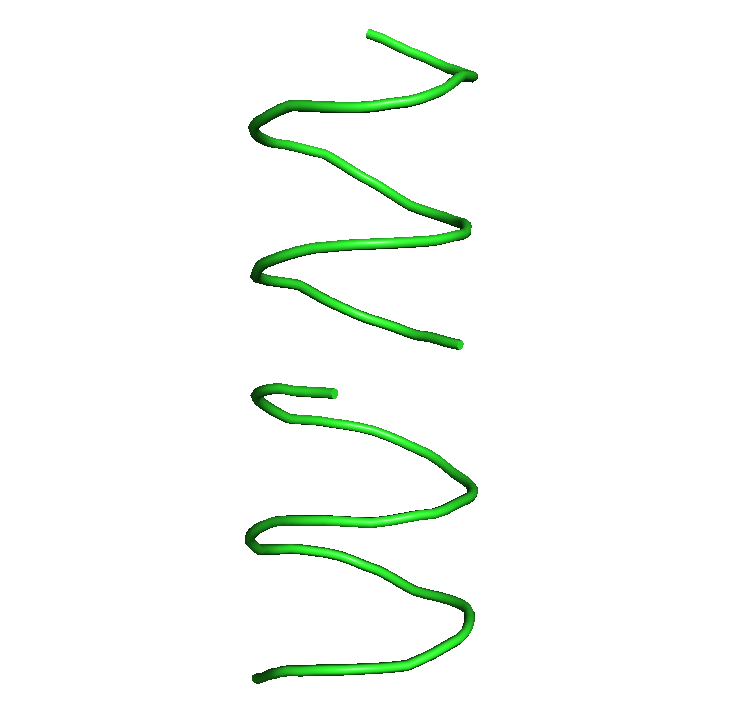}
                 {gA}
        \end{subfigure}  
                 \begin{subfigure}[b]{0.24\textwidth}
                \centering
                \includegraphics[width=\textwidth]{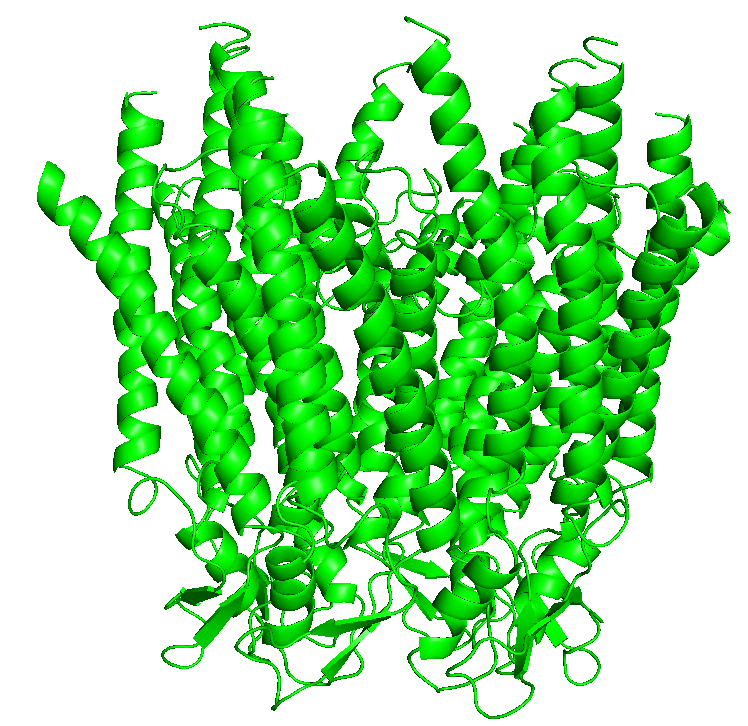}
                {Cx26}
        \end{subfigure}  
         \begin{subfigure}[b]{0.24\textwidth}
                \centering
                \includegraphics[width=\textwidth]{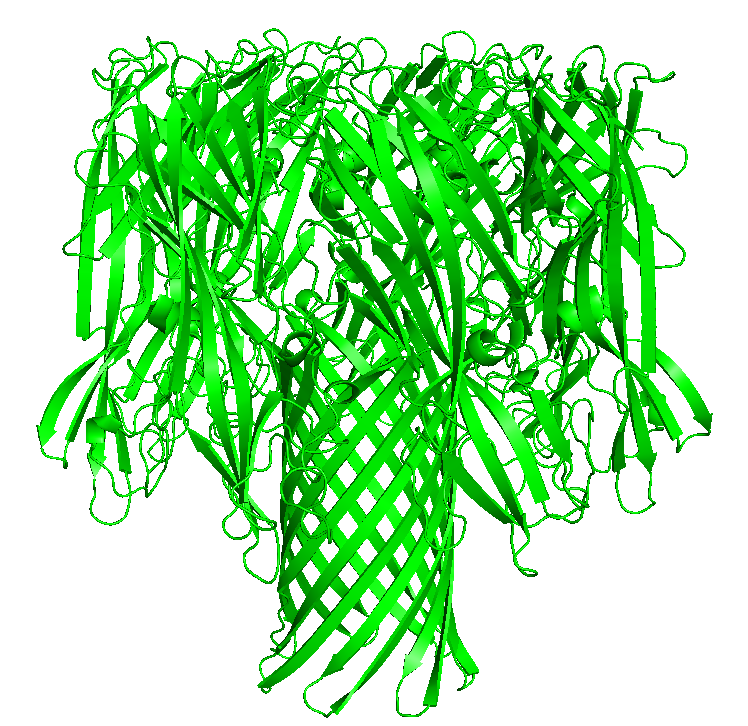}
                {$\alpha$-HL}
        \end{subfigure}  
         \begin{subfigure}[b]{0.24\textwidth}
                \centering
                \includegraphics[width=\textwidth]{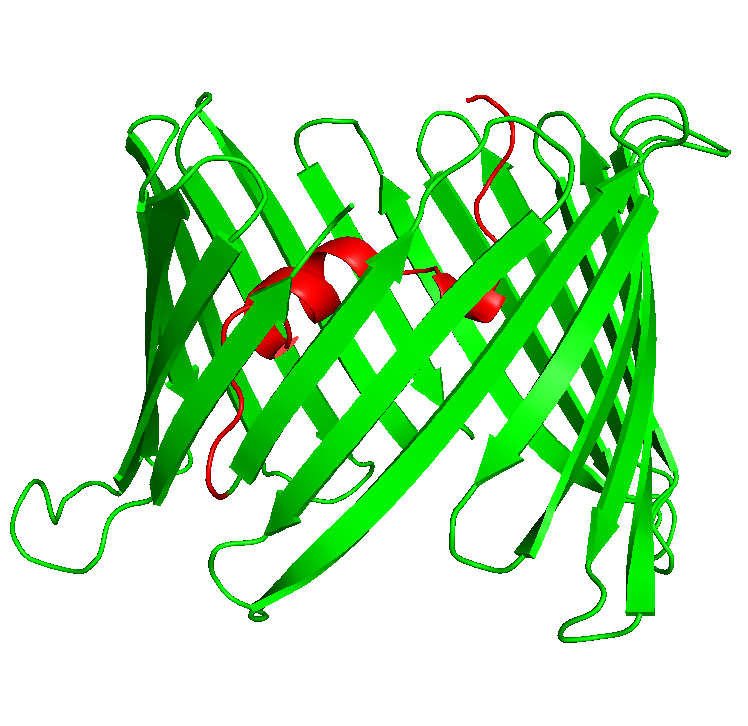}
                 {VDAC}
                  \end{subfigure}  
         \end{minipage}       
               {(a) Side view of four ion channel proteins}     
          \begin{minipage}[b]{\textwidth}
                  \begin{subfigure}[b]{0.24\textwidth}
                \centering
                \includegraphics[width=\textwidth]{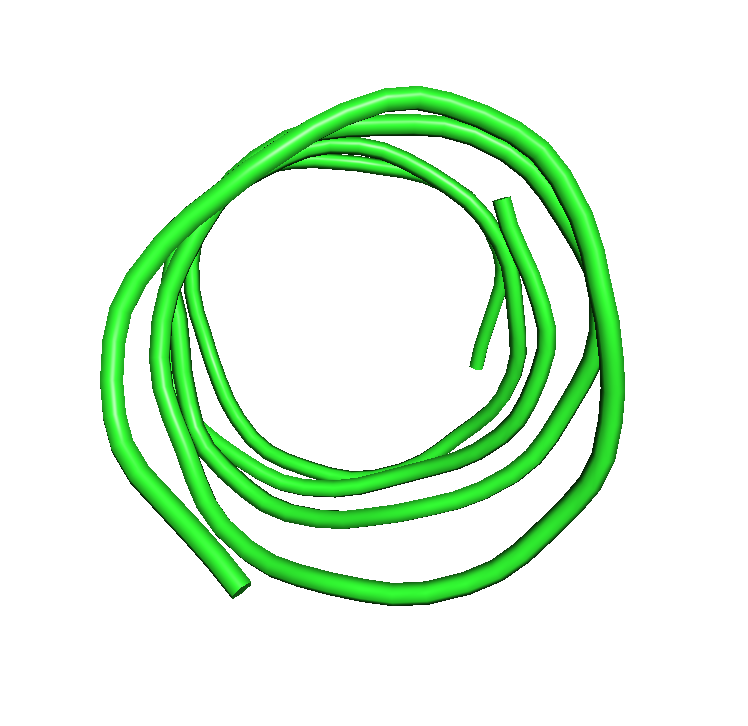}
                 {gA}
        \end{subfigure}  
                 \begin{subfigure}[b]{0.24\textwidth}
                \centering
                \includegraphics[width=\textwidth]{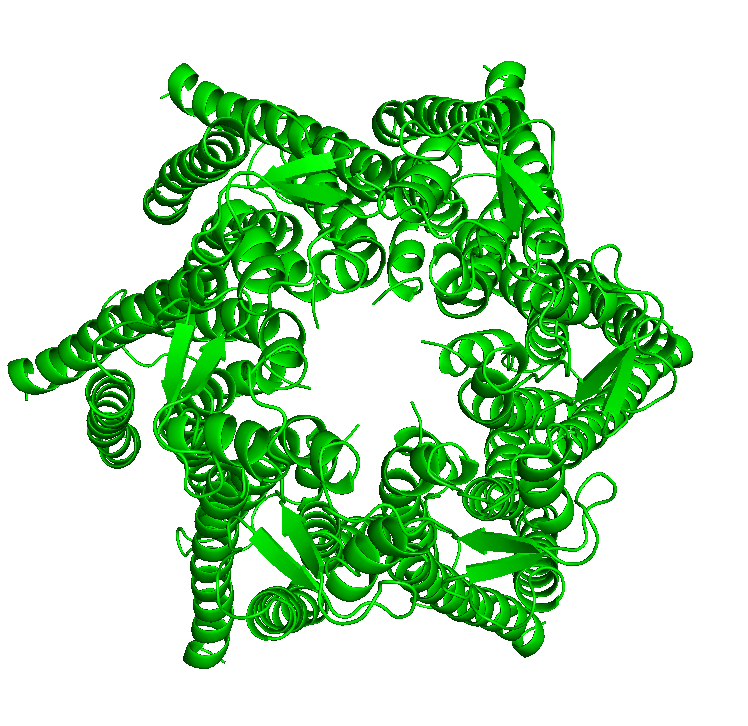}
                {Cx26}
        \end{subfigure}  
         \begin{subfigure}[b]{0.24\textwidth}
                \centering
                \includegraphics[width=\textwidth]{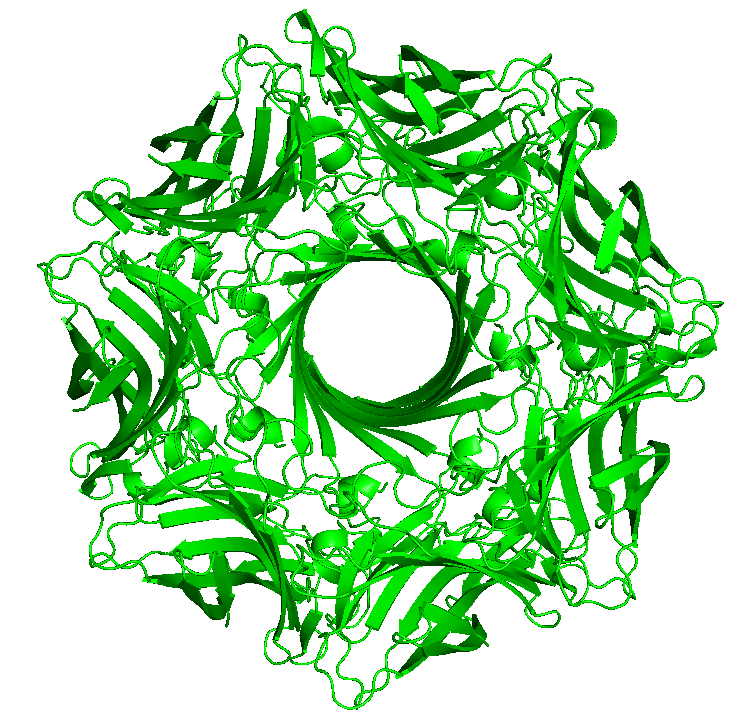}
                {$\alpha$-HL}
        \end{subfigure}  
         \begin{subfigure}[b]{0.24\textwidth}
                \centering
                \includegraphics[width=\textwidth]{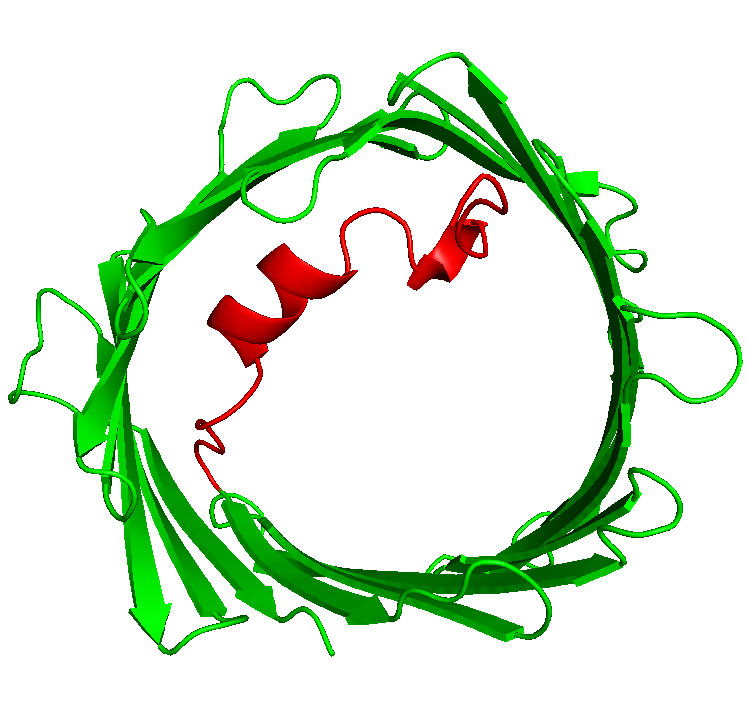}
                 {VDAC}
                  \end{subfigure} 
             \end{minipage}
                            {(b) Top view of four ion channel proteins} 
\caption{Molecular structures of the four ion channel proteins (gA, Cx26, $\alpha$-HL, and VDAC) to be used for numerical tests in Section~6  in cartoon   backbone representation --- one common way to represent a three-dimensional protein secondary structure (e.g.,  $\alpha$-helices in flat helical sheets and $\beta$-sheets in flat level sheets).} 
\label{Protein_structures}
\end{figure}

{\bf Remark:} {\em Another way to extract $D_{m,h}$ from $\hat{D}_{s,h}$ is to use  ${\cal B}_{h}$  since  ${\cal B}_{h}$ contains $D_{m,h}$ too. The reason why we use ${D}_{ms,h}$, instead of ${\cal B}_{h}$,  is to further reduce the computational cost since ${D}_{ms,h}$ contains a much smaller number of tetrahedra than ${\cal B}_{h}$.}

\section{Numerical results}

We implemented the three new schemes of Sections~3, 4, and 5  in Python based on our recent mesh work done in  \cite{xie2020finitezhen,chao2021improved} and using some mesh functions from the FEniCS project\footnote{\em  https://fenicsproject.org} and Trimesh. We then used them to modify ICMPv1 as ICMPv2. Note that ICMPv2 retains the part of ICMPv1 in the generation of  an ion channel protein molecular surface mesh $\partial D_{p,h}$ (i.e., doing so via the TMSmesh software packages developed in Lu's research group \cite{chen2012triangulated}).  Hence, both ICMPv2 and ICMPv2 are expected to generate the same ion channel protein mesh $D_{p,h}$ and expanded solvent mesh $\hat{D}_{s,h}$ when they use the same TMSmesh and TetGen parameters and the same boundary surface mesh $\partial \Omega$. Their differences mainly occur in a process of extracting membrane and solvent meshes $D_{m,h}$ and $D_{s,h}$ from an expanded solvent mesh of $\hat{D}_{s,h}$. Therefore, in this section, we mainly report the numerical results related to this extraction process.

\begin{table}[h]
\centering
\scalebox{1}{
  \begin{tabular}{|c|c|c|c|c|c|c|c|}
   \hline
  \multirowcell{2}{Ion channel\\protein}& \multirowcell{2}{Dimensions of box domain  $\Omega$ \\ $[L_{x_1}$, $L_{x_2}$; $L_{y_1}$, $L_{y_2}$; $L_{z_1}$, $L_{z_2}]$}&  \multirowcell{2}{$Z1$} &\multirowcell{2}{$Z2$} & \multirowcell{2}{$h_m$}& 
  \multirowcell{2}{$d$} & \multirowcell{2}{$c$}& \multirowcell{2}{$e$}  \\ 
  &   & &    && & &   \\ \hline
 gA &  $[-31, 31;  -30, 29; -33, 33]$ &-11     &11   &1.1  &0.5 &0.9 &0.9  \\ 
       \hline
 Cx26&   $[-67, 67;  -63, 63; -60, 62]$&-16   &16  &1.6   &0.2 &0.2 &0.8   \\ 
    \hline
    $\alpha$-HL&   $[-71, 71;  -71, 68; -39, 104]$&-11  &11  &1.1   &0.2 &0.9 &0.9   \\ 
    \hline
  VDAC&  $ [-46, 46; -46, 46;  -60, 59]$  &-12      &12 &1.2  &0.2 &0.5 &0.75 \\ 
        \hline
     \end{tabular}}
     \vspace{3mm}
  \caption{Values of box domain dimensions and main mesh parameters used in our numerical tests. Here, $d, c,$, and $e$ are three TMSmesh parameters --- $d$ is the decay rate in the Gaussian surface, $c$ is isovalue in the Gaussian surface, and $e$ is an approximation precision between trilinear surface and Gaussian surface.}
   \label{table:parameters4NumericalTests}
\end{table}

In particular, we did numerical tests on four ion channel proteins: (1) A gramicidin A (GA), (2) a connexin 26 gap junction channel (Cx26), (3) a  staphylococcal $\alpha$-hemolysin ($\alpha$-HL), and (4) a  voltage-dependent anion channel (VDAC). Their crystallographic three-dimensional molecular structures  can be downloaded from the  Protein Data Bank\footnote{\em https://www.rcsb.org} with the PDB identification numbers 1GRM,  2ZW3, 7AHL, and 3EMN, respectively.  But, in this work, we  downloaded them from the Orientations of Proteins in Membranes (OPM) database\footnote{\em https://opm.phar.umich.edu} since these molecular structures  have satisfied our assumptions made in the partition \eqref{DomainPartition}. That is, the protein structure has been transformed such that the normal direction of the top membrane surface is in the $z$-axis direction and the membrane location numbers $Z1$ and $Z2$ are given. See Table~\ref{table:parameters4NumericalTests} for their values.  

Figure~\ref{Protein_structures} displays these four ion channel proteins in cartoon  views. The $\alpha$-helix of VDAC has been colored in red to more clearly view it  in Figure~\ref{Protein_structures}. From these plots we can see that the ion channel pores of these four  proteins have different geometrical complexities. Thus, these four proteins are good for numerical tests on the efficiency of our new numerical schemes and a comparison study between ICMPv1 and ICMPv2.
                        
Table~\ref{table:parameters4NumericalTests} lists the values of box domain dimensions, three parameters $Z1, Z2$, and $h_m$ from our new schemes, and three parameters $h, d$, and $c$  from the molecular surface software TMSmesh. In the numerical tests,  we fixed the other parameter values of our schemes as follows:
 \[    \eta_1 = 20, \quad \eta_2 = 20, \quad \eta_3 = 20, \quad h_s=2h_m, \quad \tau=h_m.\]
We also fixed the command line switches of TetGen as {\em `-q1.2aVpiT1e-10AAYYCnQ'}, whose definitions and usages can be found in TetGen's manual webpage\footnote{{\em https://wias-berlin.de/software/tetgen/1.5/doc/manual/manual.pdf}}. Thus, we did not list them in Table~\ref{table:parameters4NumericalTests}. Actually, all the box domain dimensions of Table~\ref{table:parameters4NumericalTests} can be produced from the formulas of \eqref{eq:boxDomainSize}. Even so, we have listed them in Table~\ref{table:parameters4NumericalTests} for clarity. 

 \begin{table}[h]
\centering
  \begin{tabular}{|c|c|c|c|c|c|c|}
   \hline
  \multirowcell{2}{Channel\\ protein}&\multicolumn{6}{c|}{Number of vertices }        \\  
  \cline{2-7}
  &   $\hat{D}_{s,h}$     &$D_{s,h}^p$&$\Omega_h$ &  $D_{s,h}$  &  $D_{m,h}$&  $D_{p,h}$ \\ \hline
   \multicolumn{7}{|c|}{Mesh data generated by ICMPv1}\\
  \hline
gA &17711&662  &24155   &9359  &10653 &11012     \\ 
    \hline
   Cx26  &74185&1091 &156182   &55869  &24304 &106354    \\ 
       \hline
   $\alpha$-HL  &109911&1615 &230657   &  101212&13582 &160141   \\ 
    \hline
VDAC   &32842 &2323&52640    &23877  &12981 &28863    \\ 
    \hline
    \multicolumn{7}{|c|}{Mesh data generated by ICMPv2}\\
 \hline
gA & 27075&690  &33612    & 15254  &15488 & 11105   \\ 
    \hline
   Cx26  &73984&1096 &   158396&  55255& 22645&  108769   \\ 
       \hline
   $\alpha$-HL  &116649&1545 &237535   &105027  &16757 &  160281 \\ 
    \hline
VDAC   &35681&2426  &56584     &25538  &13440 &29968   \\ 
    \hline
\end{tabular}

\vspace{3mm}

  \begin{tabular}{|c|c|c|c|c|c|c|}
   \hline
  \multirowcell{2}{Channel\\ protein}&\multicolumn{6}{c|}{Number of tetrahedra }        \\  
  \cline{2-7}
  &   $\hat{D}_{s,h}$     &$D_{s,h}^p$&$\Omega_h$ &  $D_{s,h}$  &  $D_{m,h}$&  $D_{p,h}$ \\ \hline
   \multicolumn{7}{|c|}{Mesh data generated by ICMPv1}\\
  \hline     
gA  & 95145& 2488&149558 &44092  &51053 &54413   \\ 
    \hline
   Cx26    &387042 &4221  &979905 & 274291 & 112751 &592863 \\ 
       \hline
   $\alpha$-HL  &563949 &6471  &1448649 & 501402 & 62547 &884700 \\ 
    \hline
VDAC   &176335   &9788&328992 &116287  &60048 &152657   \\ 
    \hline
    \multicolumn{7}{|c|}{Mesh data generated by ICMPv2}\\
 \hline
gA &143799  &2579&198776 &70787  &73012 &54977   \\ 
    \hline
   Cx26   &375132 &4282  &983421 &  270715&104417 &608289   \\ 
       \hline
   $\alpha$-HL  &597027 &6403  &1482295 & 519312 &77715  &885268 \\ 
    \hline
VDAC   &183414   &10334&343072 &159658  &60685 &159658    \\ 
    \hline
\end{tabular}

\vspace{3mm}

  \caption{A comparison of mesh data generated by ICMPv1 with those by ICMPv2. }
   \label{table:meshData}
\end{table}

\begin{figure}[h!]
        \centering
                \begin{minipage}[b]{\textwidth}
                 \begin{subfigure}[b]{0.24\textwidth}
                \centering
                \includegraphics[width=\textwidth]{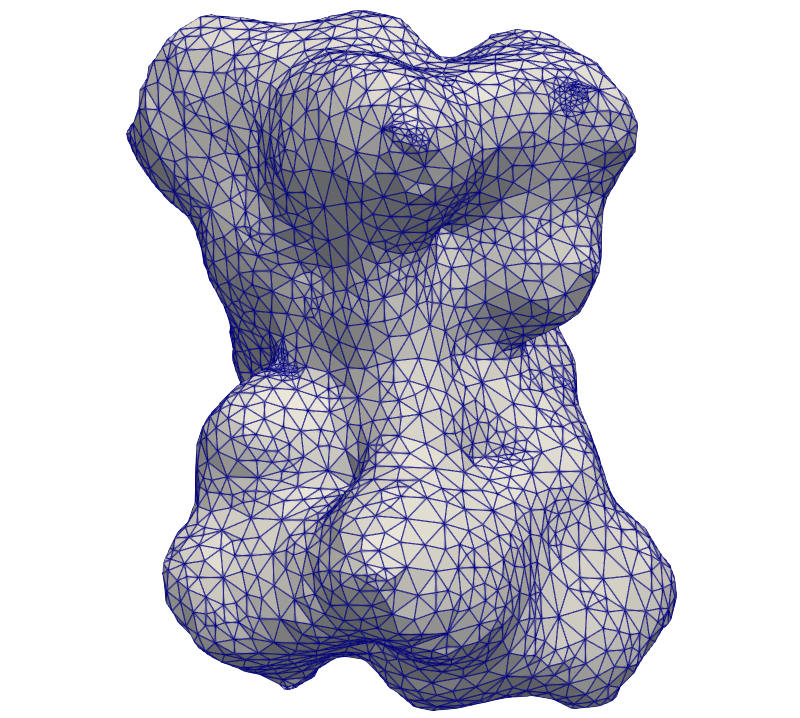}
                 {gA}
              \label{fig:1GRM}
        \end{subfigure}  
                 \begin{subfigure}[b]{0.24\textwidth}
                \centering
                \includegraphics[width=\textwidth]{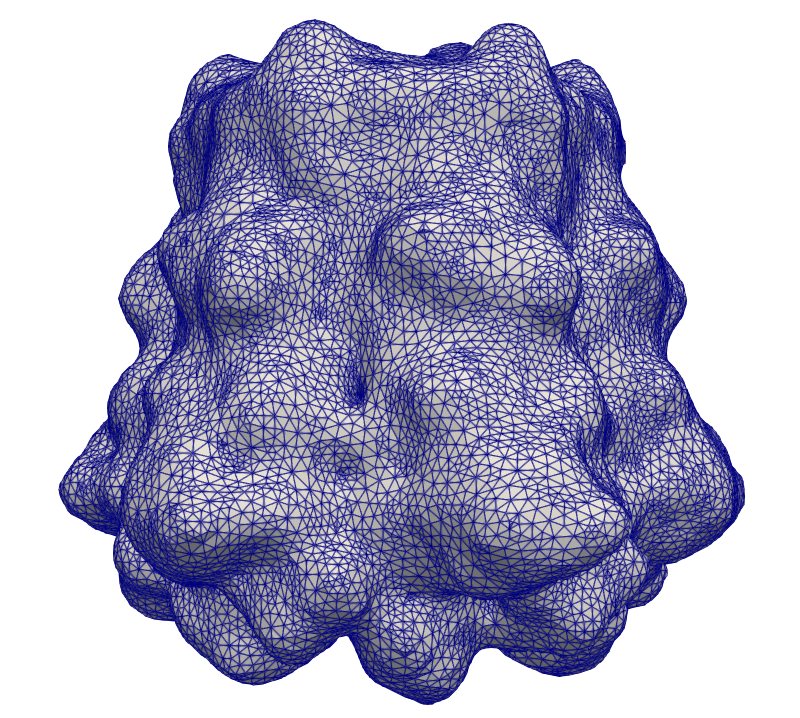}
                {Cx26}
             \label{fig:2zw3HalDpMeshf}
        \end{subfigure}  
         \begin{subfigure}[b]{0.24\textwidth}
                \centering
                \includegraphics[width=\textwidth]{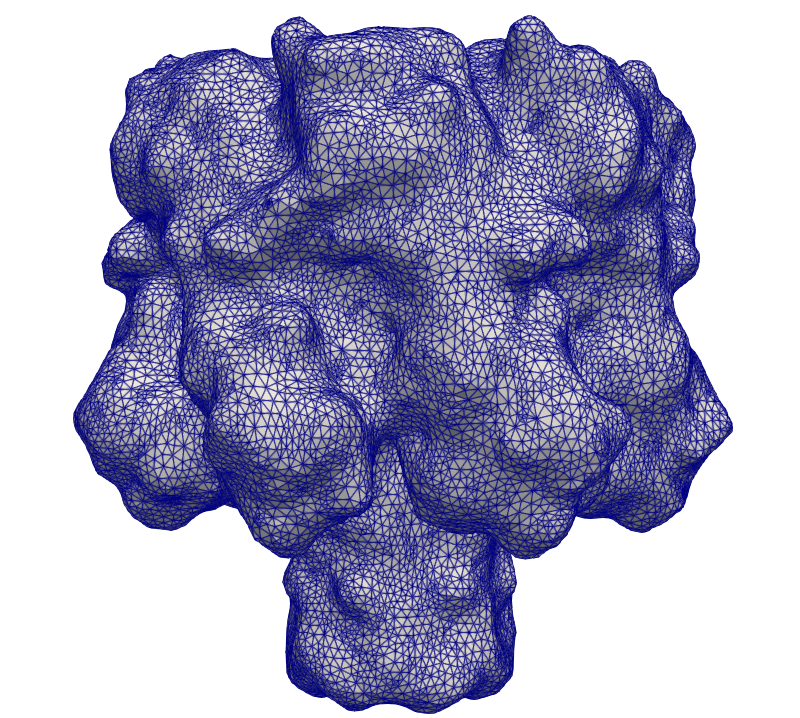}
                {$\alpha$-HL}
              \label{fig:7AHL}
        \end{subfigure}  
         \begin{subfigure}[b]{0.24\textwidth}
                \centering
                \includegraphics[width=\textwidth]{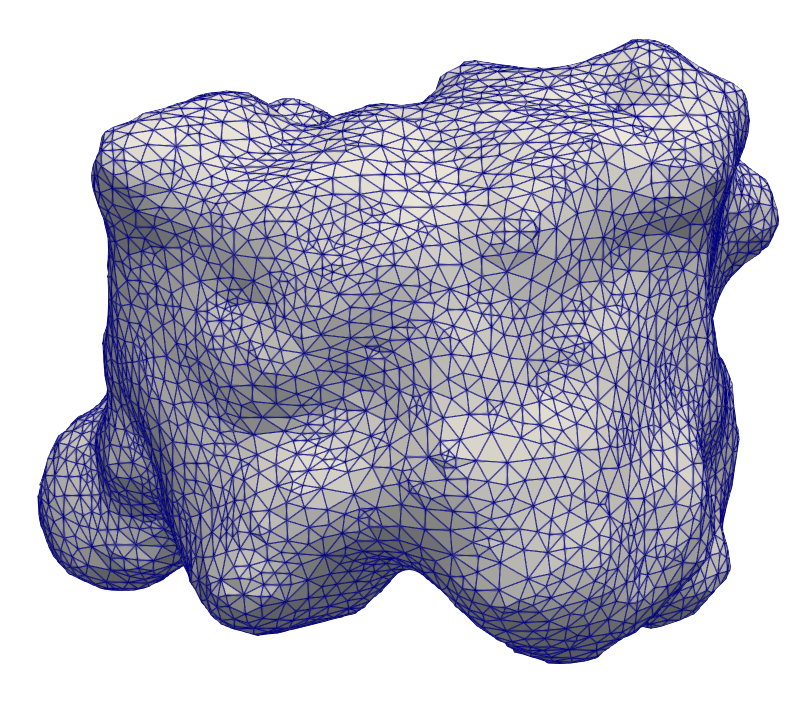}
                 {VDAC}
              \label{fig:3emne_2}
                      \end{subfigure} 
                      \end{minipage}
                      {(a) Side view of four ion channel protein meshes}    
                \begin{minipage}[b]{\textwidth}
                    \begin{subfigure}[b]{0.24\textwidth}
                \centering
                \includegraphics[width=\textwidth]{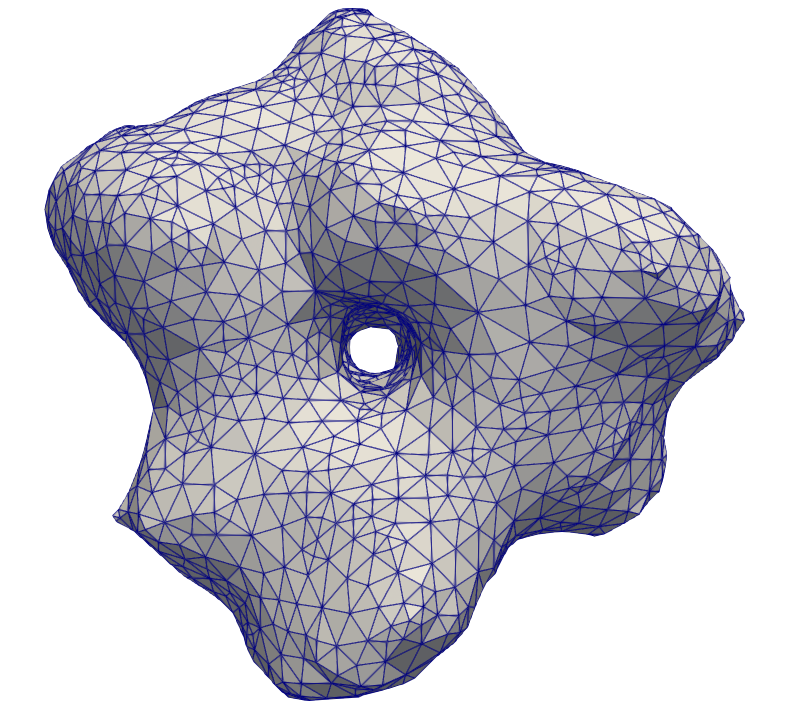}
                 {gA}
              \label{fig:1GRM_2}
        \end{subfigure}  
                 \begin{subfigure}[b]{0.24\textwidth}
                \centering
                \includegraphics[width=\textwidth]{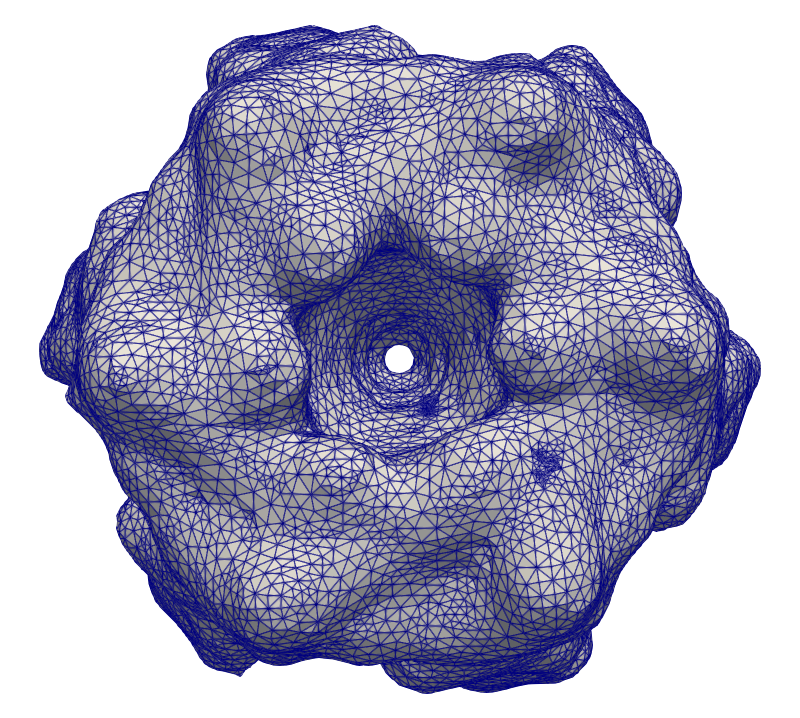}
                {Cx26}
             \label{fig:2zw3HalDpMeshf_2}
        \end{subfigure}  
         \begin{subfigure}[b]{0.24\textwidth}
                \centering
                \includegraphics[width=\textwidth]{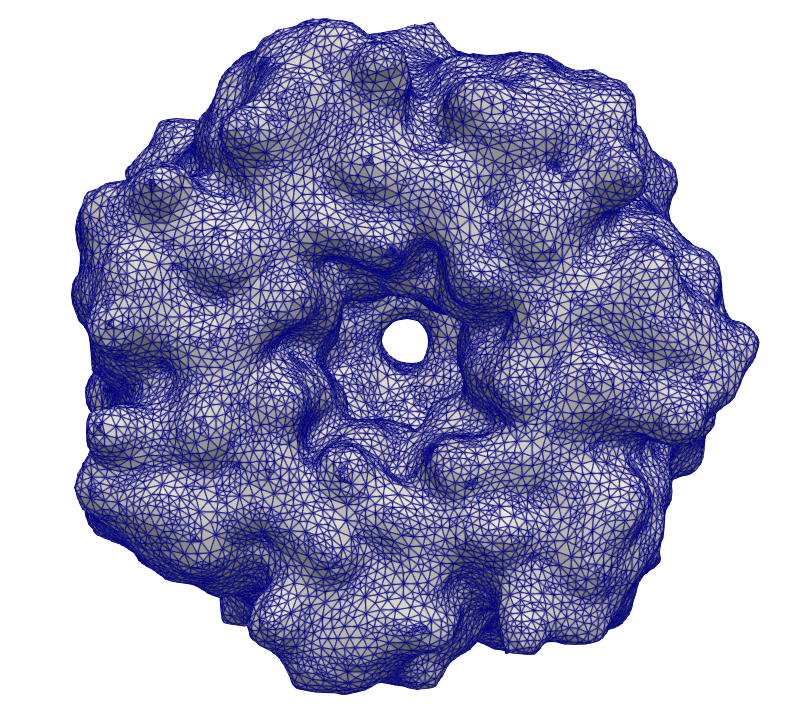}
                {$\alpha$-HL}
              \label{fig:7AHL_a}
        \end{subfigure}  
         \begin{subfigure}[b]{0.24\textwidth}
                \centering
                \includegraphics[width=\textwidth]{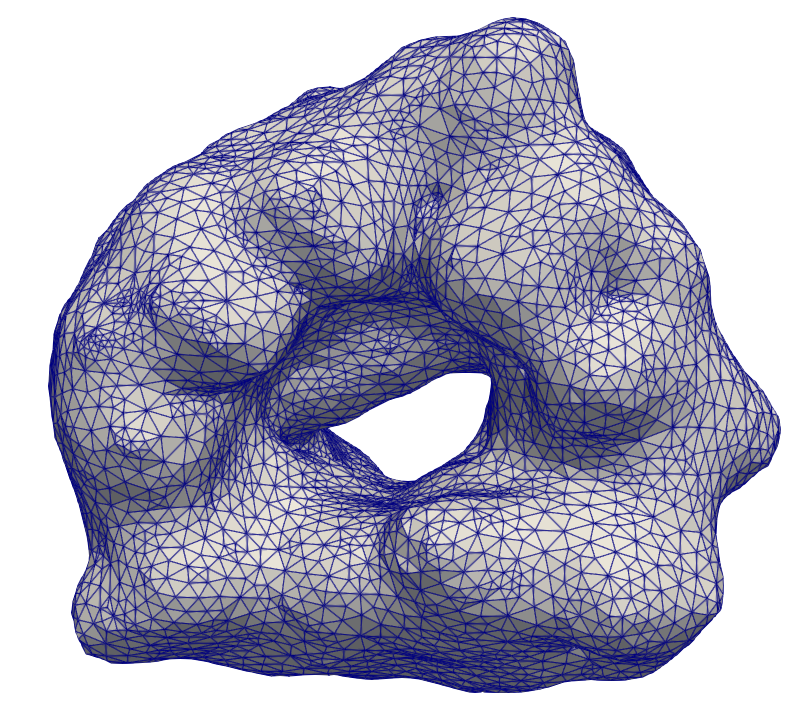}
                 {VDAC}
              \label{fig:3emne_3}
                      \end{subfigure}  
               \end{minipage}   
               {(b) Top view of four ion channel protein meshes}    
\caption{Top views of the four ion channel protein meshes $D_{p,h}$ generated by ICMPv2 according to the four molecular structures displayed in Figure~\ref{Protein_structures}.}
\label{proteinMeshes}
\end{figure}

\begin{figure}[pt!]
        \centering
                 \begin{subfigure}[b]{0.4\textwidth}
                \centering
               \includegraphics[width=\textwidth]{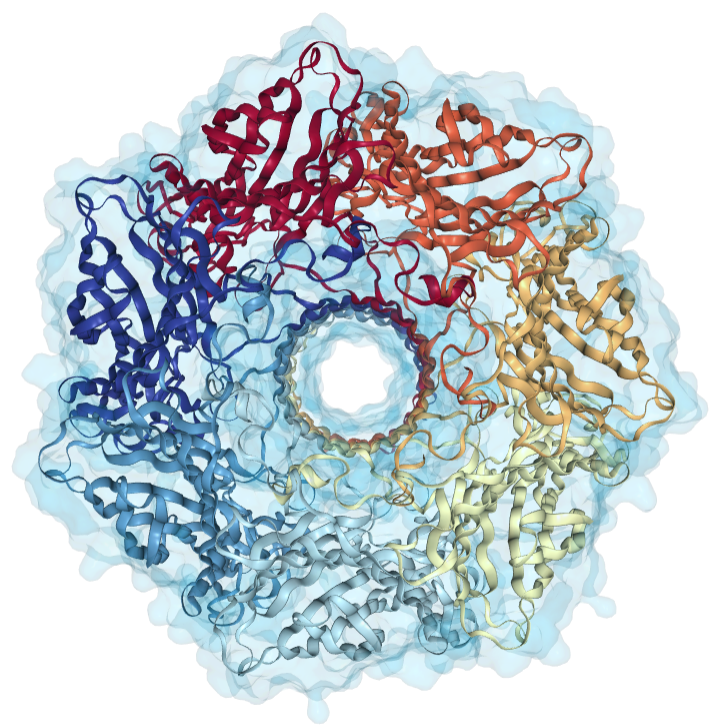}
               \caption{A top view  }
               \label{fig:7AHL_view1}
        \end{subfigure}  
         \qquad \qquad
         \begin{subfigure}[b]{0.4\textwidth}
                \centering
                \includegraphics[width=\textwidth]{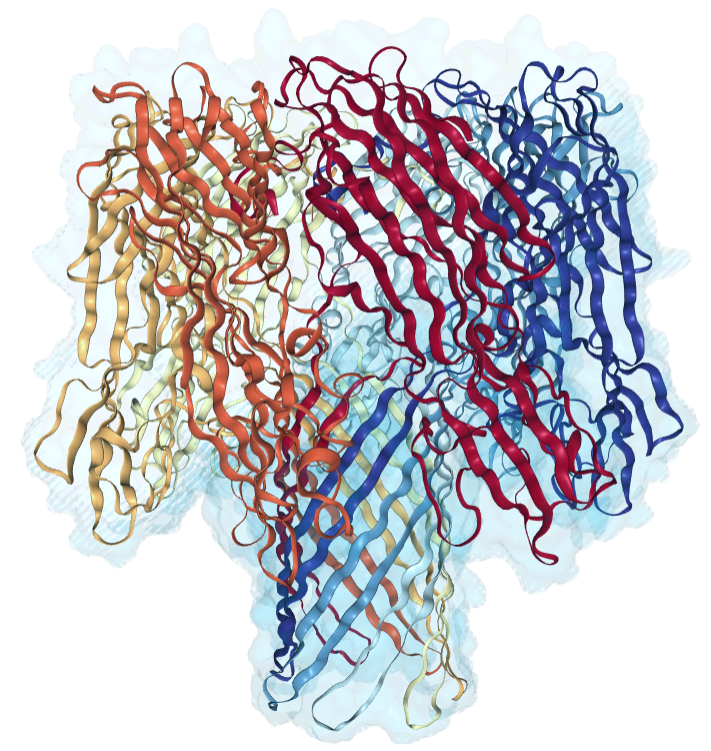}
                \caption{A side view}
              \label{fig:7AHL_view2}
        \end{subfigure}           
\caption{A comparison of  a protein mesh ${D}_{p,h}$ generated by ICMPv2 
with a molecular structure of $\alpha$-HL (PDB ID: 7AHL). Here the protein structure is depicted in cartoon   backbone representation and colored in red, blue, yellow, and cyan.}    
\label{fig:1bl8meshMolecule}
\end{figure}

\begin{figure}[pht!]
        \centering
         \begin{subfigure}[b]{0.45\textwidth}
                \centering
               \includegraphics[width=\textwidth]{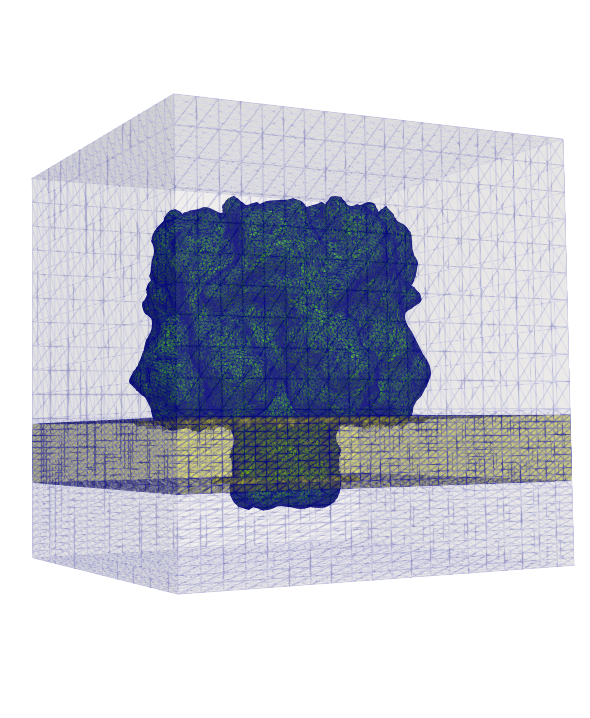}
                \caption{A view of $\Omega_h$ }
             \label{fig:1bl8DsSide-17}
        \end{subfigure} 
        \qquad 
        \begin{subfigure}[b]{0.45\textwidth}
                \centering
                \includegraphics[width=\textwidth]{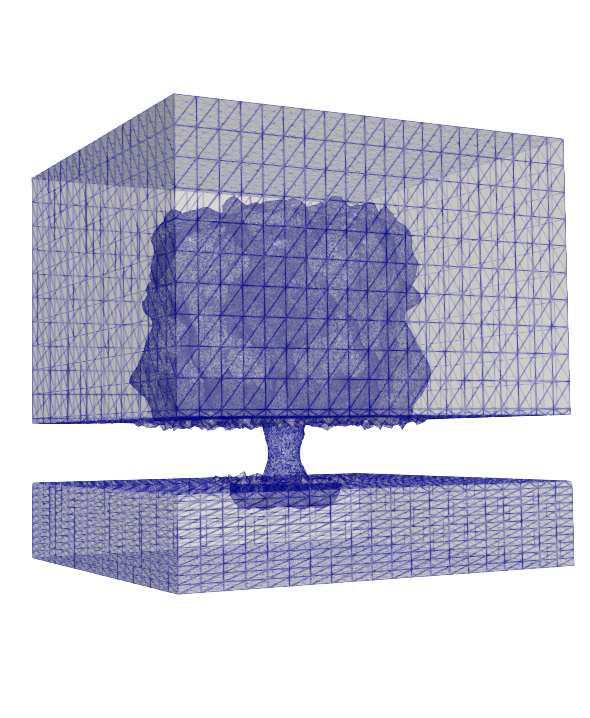}
                \caption{Side view of $D_{s,h}$ }
             \label{fig:1bl8DsSide}
        \end{subfigure} 
                  \begin{subfigure}[b]{0.45\textwidth}
                \centering
               \includegraphics[width=\textwidth]{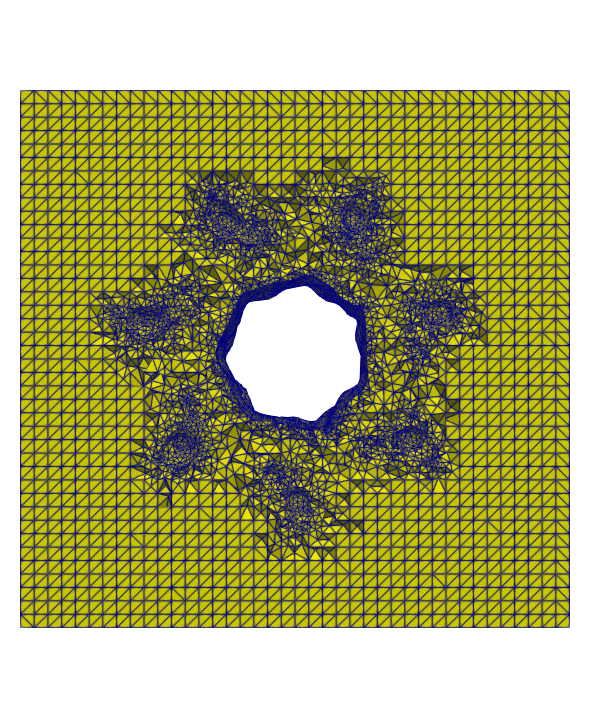}
                \caption{Top view of ${D}_{m,h}$ at $z=Z2$}
               \label{fig:3emnMembraneSubmesh2}
        \end{subfigure}  
                 \qquad
         \begin{subfigure}[b]{0.45\textwidth}
                \centering
                \includegraphics[width=0.98\textwidth]{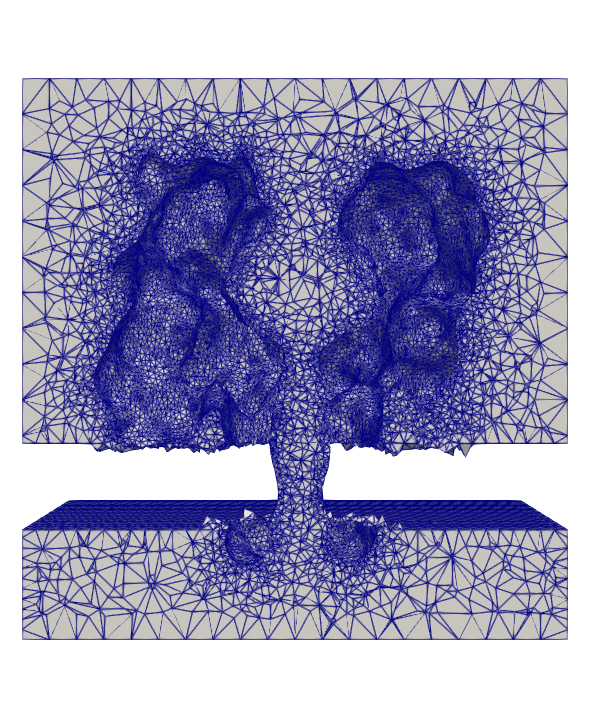}
                \caption{Cross section of $D_{s,h}$ on the $yz$-plane}
             \label{fig:1bl8DsPore}
        \end{subfigure} 
\caption{The whole domain mesh $\Omega_{h}$, solvent mesh $D_{s,h}$, and membrane mesh $D_{m,h}$ generated by ICMPv2 for the ion channel protein $\alpha$-HL. }
\label{fig:7AHL_meshes}
\end{figure}

\begin{figure}[h!]
        \centering
        A side view of a membrane mesh $D_{m, h}$
        
                 \begin{subfigure}[b]{0.4\textwidth}
                \centering
                \includegraphics[width=\textwidth]{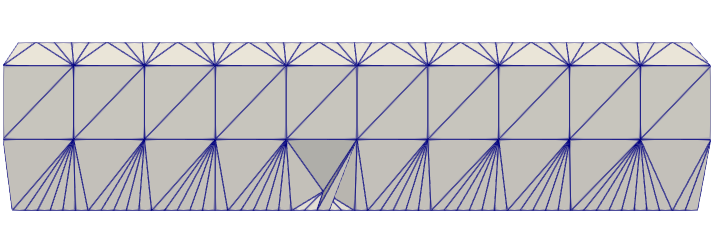}
                 \caption{By ICMPv1}
              \label{fig:3emnDmMesh_luSide}
                      \end{subfigure}    
                      \qquad       \qquad   
                   \begin{subfigure}[b]{0.4\textwidth}
                \centering
                \includegraphics[width=\textwidth]{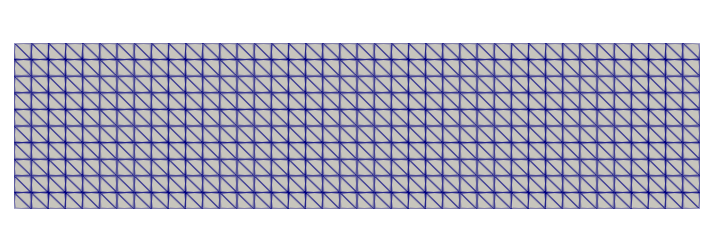}
                 \caption{By ICMPv2}
              \label{fig:3emnDmMesh_zhenSide}
        \end{subfigure}  
        
        A top view of a membrane mesh $D_{m, h}$
        
         \begin{subfigure}[b]{0.4\textwidth}
                \centering
                \includegraphics[width=0.9\textwidth]{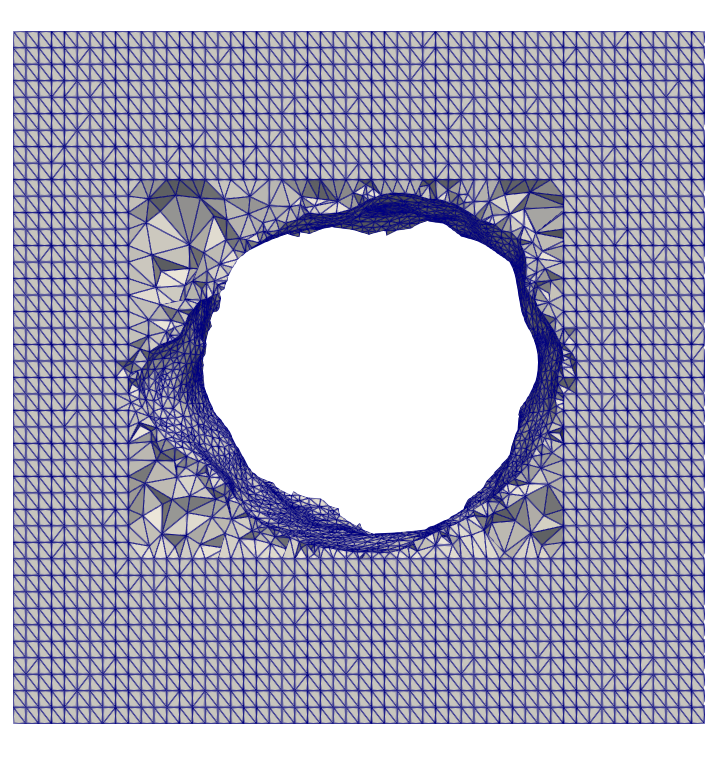}
                \caption{By ICMPv1}
              \label{fig:3emnDmMesh_luTop}
        \end{subfigure}  
        \qquad. \qquad
                      \begin{subfigure}[b]{0.4\textwidth}
                \centering
                \includegraphics[width=0.9\textwidth]{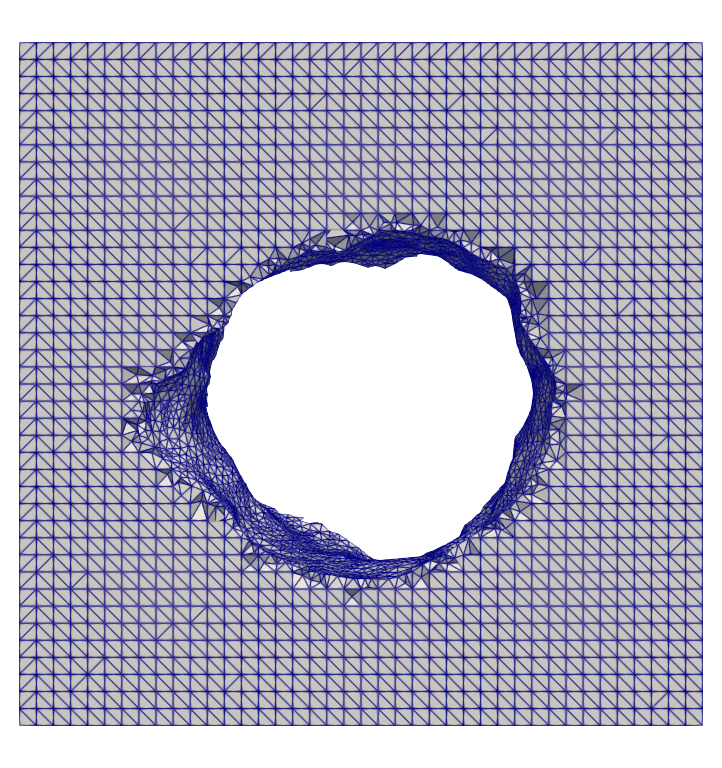}
                 \caption{By ICMPv2}
              \label{fig:3emnDmMesh_zhenTop}
                      \end{subfigure}   
                      
           A view of the bottom portion $D_{s, h}^b$ of a solvent mesh $D_{s,h}$
           
         \begin{subfigure}[b]{0.4\textwidth}
                \centering
                \includegraphics[width=\textwidth]{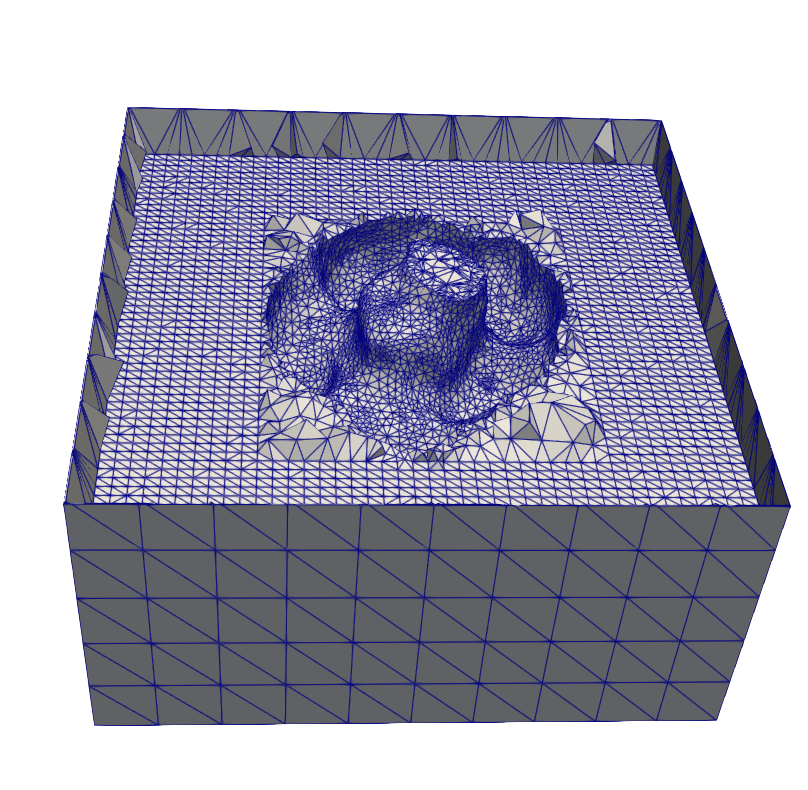}
                 \caption{By ICMPv1}
              \label{fig:3emnDsMesh_luSide}
        \end{subfigure} 
        \qquad \qquad
                 \begin{subfigure}[b]{0.4\textwidth}
                \centering
                \includegraphics[width=\textwidth]{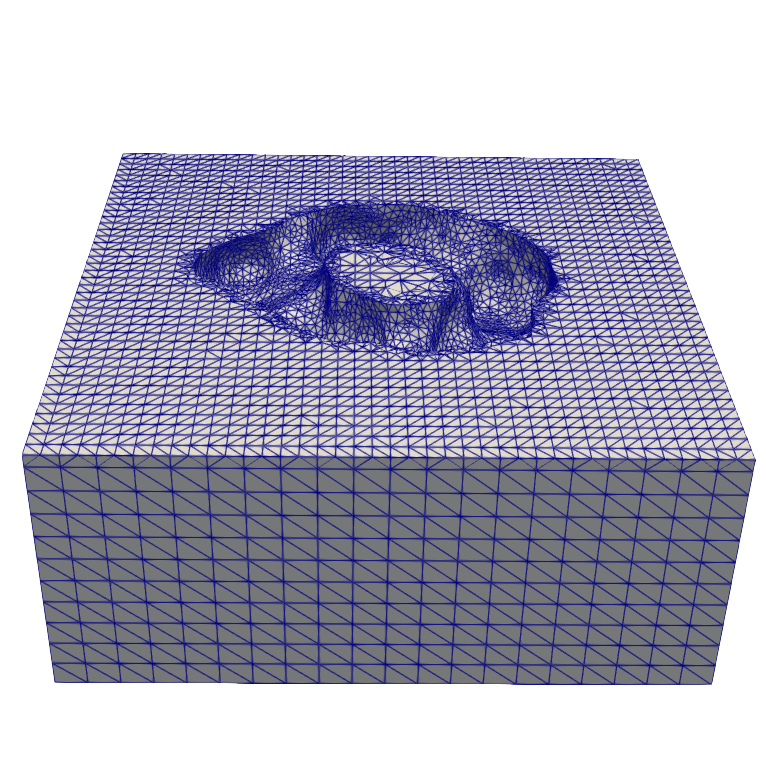}
                \caption{By ICMPv2}
             \label{fig:3emnDsMesh_zhenSide}
        \end{subfigure} 
\caption{A comparison of the membrane and solvent meshes $D_{m,h}$ and $D_{s,h}$ generated by ICMPv1 with those by ICMPv2 for VDAC. Here $D_{s,h}^b$ is defined in \eqref{solvenDomainPartition}.}
\label{fig:3emnMeshProfileComparison}
\end{figure}

We did all the numerical tests on a MacBook Pro computer with one 2.6 GHz Intel core i7 processor and 16 GB memory. We listed the mesh data generated from  ICMPv1 and ICMPv2  in Table~\ref{table:meshData} and reported other test results in Table~\ref{table:meshDataCPU} and Figures~\ref{proteinMeshes} to \ref{fig:3emnMeshProfileComparison}.

Figure~\ref{proteinMeshes} displays the triangular surface meshes of the four ion channel proteins generated from ICMPv2 through using the software package TMSmesh. Note that VDAC has a much larger channel pore than the others. 

Because $\alpha$-HL has a complex molecular structure, we take it as an example to show how an ion channel protein mesh $D_{p,h}$ generated by ICMPv2 to fit a molecular structure.  From 
Figure~\ref{fig:1bl8meshMolecule} it can be seen that the protein mesh  $D_{p,h}$ can hold the molecular structure very well.
We further display the box domain mesh $\Omega_{h}$, solvent mesh $D_{s,h}$, and membrane mesh $D_{m,h}$ as well as a cross section view of solvent mesh $D_{s,h}$ in Figure~\ref{fig:7AHL_meshes}. From this figure it can be seen that the unstructured tetrahedral meshes generated from ICMPv2 can  catch well the main features of the complicated geometric shapes of the protein, membrane, and solvent regions $D_p, D_m$, and $D_s$ and the complex interfaces of the box domain mesh $\Omega_{h}$. 

\begin{figure}[htp!]
        \centering
         \begin{subfigure}[b]{0.34\textwidth}
                \centering
                \includegraphics[width=\textwidth]{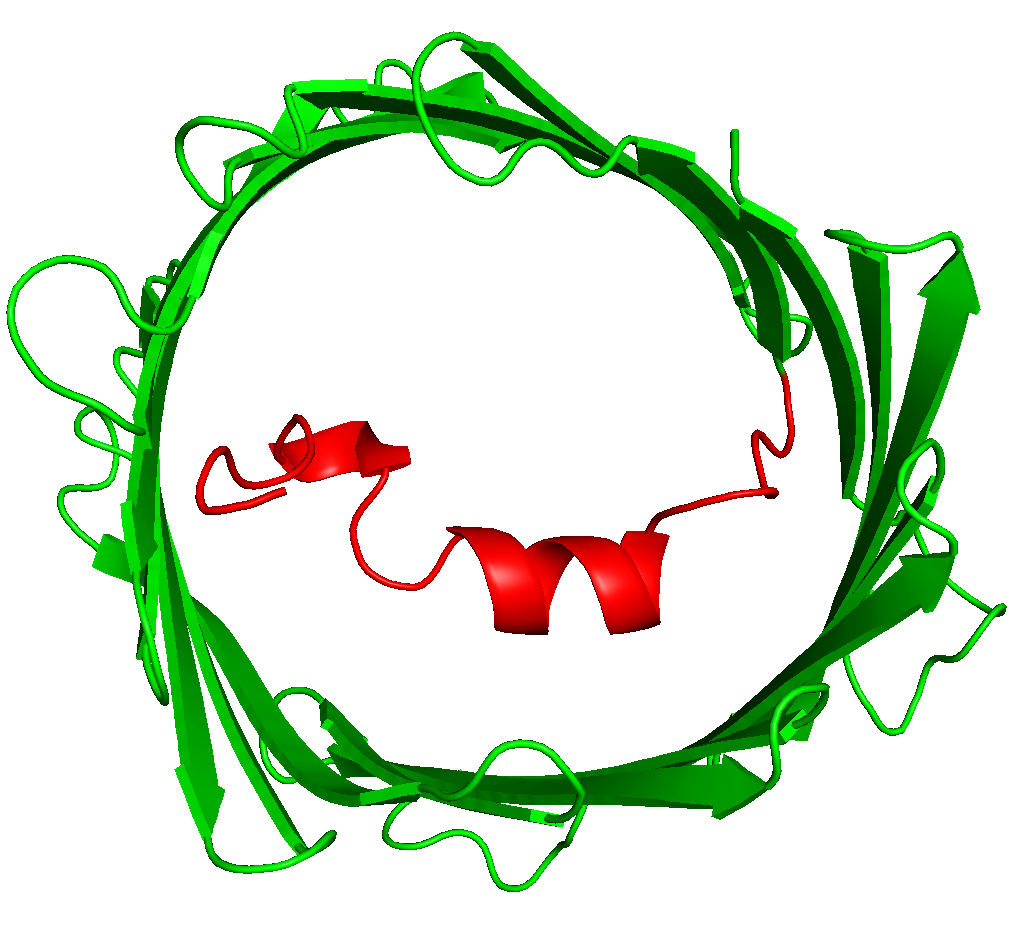}
                \caption{VDAC molecular structure}
              \label{fig:m3emnStructure}
        \end{subfigure}  
        \qquad
         \begin{subfigure}[b]{0.34\textwidth}
                \centering
                \includegraphics[width=\textwidth]{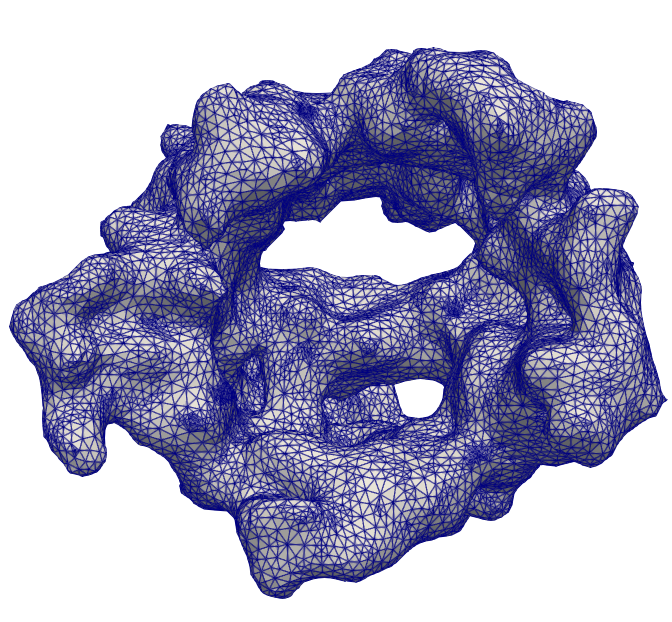}
                 \caption{Protein mesh $ D_{p,h}$}
              \label{fig:m3emnMesh}
        \end{subfigure}   
        
\vspace{3mm}

        \centering
         A top view of a  membrane mesh $D_{m, h}$
         
              \begin{subfigure}[b]{0.34\textwidth}
                \centering
                \includegraphics[width=\textwidth]{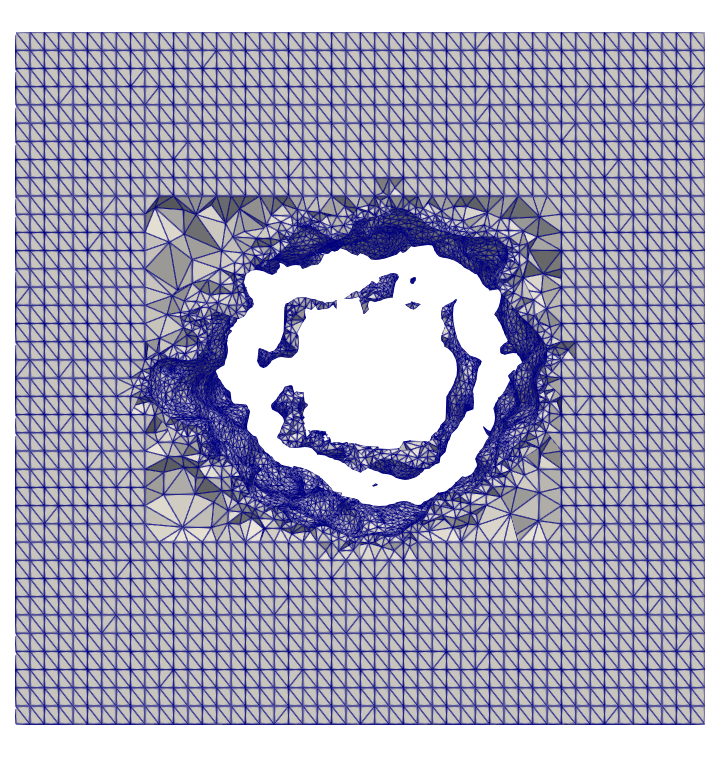}
                 \caption{By ICMPv1 }
             \label{fig:3emnRotateZ20DmMesh_luTop}
                      \end{subfigure}   
        \quad \quad
                      \begin{subfigure}[b]{0.34\textwidth}
                \centering
                \includegraphics[width=\textwidth]{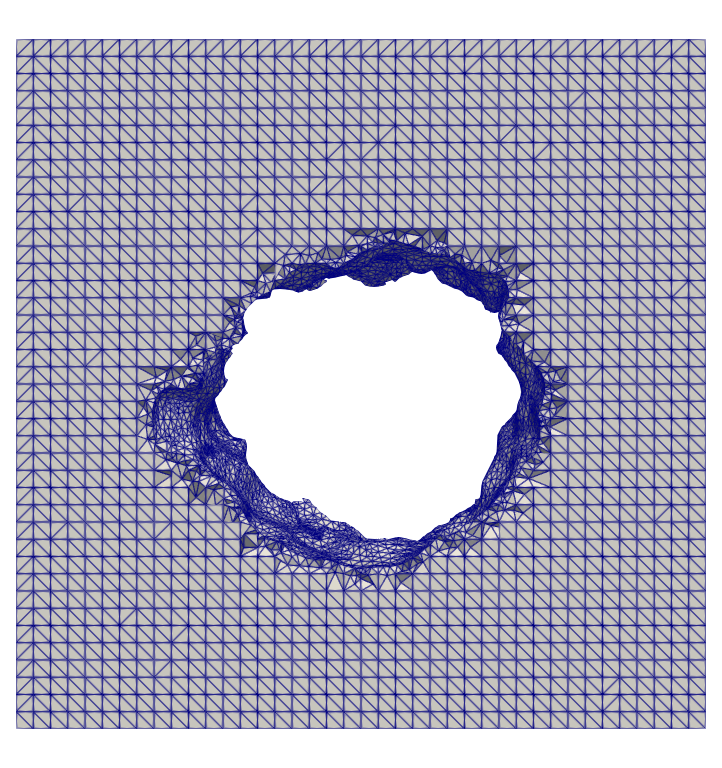}
                 \caption{By ICMPv2}
             \label{fig:3emnRotateZ20DmMesh_zhenTop}
                      \end{subfigure}   
                      
                A side view of a portion  $D_{s,h}^p$ of $D_{s,h}$ within the  ion channel pore
                
        \begin{subfigure}[b]{0.34\textwidth}
                \centering
                \includegraphics[width=\textwidth]{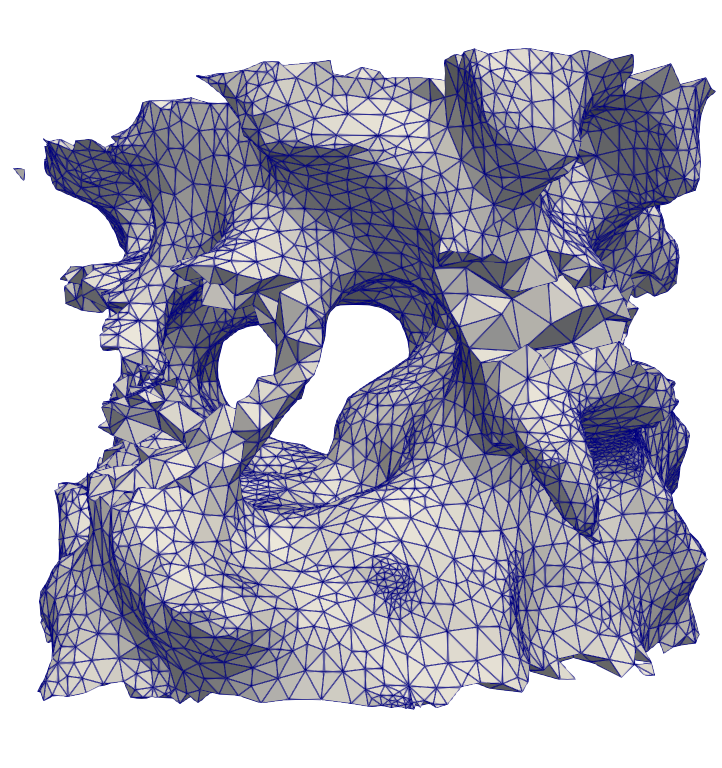}
                \caption{By ICMPv1 }
             \label{fig:3emnRotateZ20ChannelPoreMesh_lu}
        \end{subfigure} 
         \quad \quad
                 \begin{subfigure}[b]{0.34\textwidth}
                \centering
                \includegraphics[width=\textwidth]{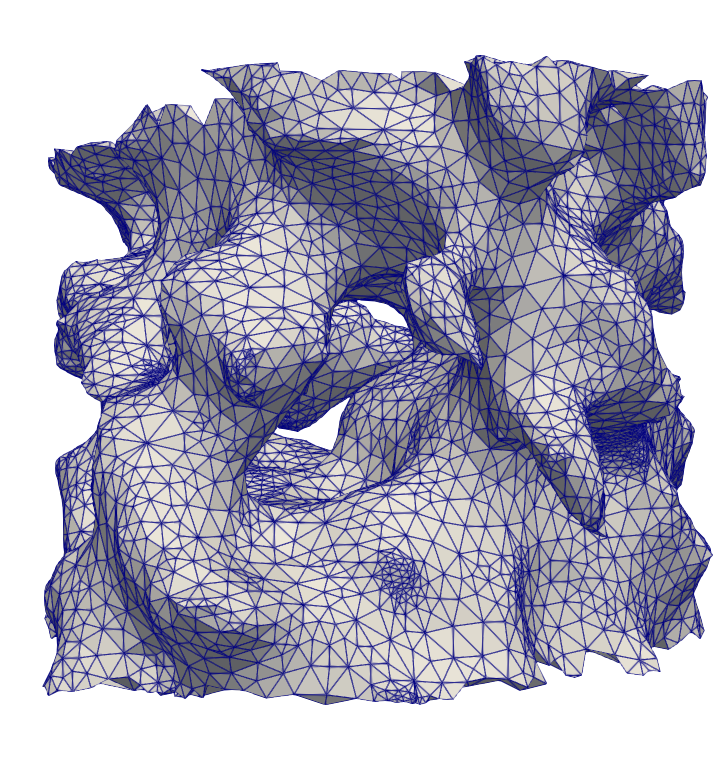}
                \caption{By ICMPv2}
             \label{fig:3emnRotateZ20ChannelPoreMesh_zhen2}
        \end{subfigure} 
\caption{(a,b): A molecular structure of mVDAC and a protein mesh of mVDAC generated by ICMPv2. (c - f): A comparison of the tetrahedral meshes $D_{m,h}$ and $D_{s,h}^p$ generated by ICMPv2 with those by ICMPv1 for  mVDAC. Here  $D_{s,h}^p$ is defined in \eqref{solvenDomainPartition}.}
\label{fig:3emnRotateZ20Comparison} 
\end{figure}

Figure~\ref{fig:3emnMeshProfileComparison} presents a comparison of the membrane and solvent meshes ${D}_{m,h}$ and ${D}_{s,h}$ generated by ICMPv1 and ICMPv2 for VDAC. From it we can see that ICMPv2 can generate either ${D}_{m,h}$ or ${D}_{s,h}$ in  higher quality than ICMPv1. We should point out that we spent a lot of time on the adjustment of mesh parameters to let ICMPv1 be able to extract ${D}_{m,h}$ and ${D}_{s,h}$ from the given expanded solvent mesh $\hat{D}_{s,h}$ successively in the sense that the membrane mesh ${D}_{m,h}$ does not contain any tetrahedron from the solvent mesh  ${D}_{s,h}$. Even so, the numerical schemes of ICMPv1 for constructing a box surface mesh and for selecting a set of mesh points from the bottom and top membrane surfaces still caused extraction problems, which decayed mesh quality.

 \begin{table}[h]
\centering
  \begin{tabular}{|c|c|c|c|c|}
   \hline
  \multirowcell{2}{Mesh\\ package}&\multicolumn{4}{c|}{CPU time (in seconds)  }   \\  
  \cline{2-5}
&  gA    & Cx26 & $\alpha$-HL & VDAC  \\ \hline
 ICMPv1  & 12.3  & 90.4 &94.5&50.1   \\ \hline
 ICMPv2  & 1.1  &2.2  &2.5& 1.5   \\ \hline
     \end{tabular}
     \vspace{3mm}
  \caption{A comparison of the computer performance of ICMPv1 with that of  ICMPv2.}
   \label{table:meshDataCPU}
\end{table}

Table~\ref{table:meshDataCPU} reports the  performance of ICMPv1 and ICMPv2 in computer CPU time spent on the extraction of membrane and solvent meshes $D_{p,h}$ and $D_{s,h}$ from a given expanded solvent mesh $\hat{D}_{s,h}$. From Table~\ref{table:meshDataCPU} it can be seen that our new numerical algorithms reported in Sections 3, 4, and 5 can be much more efficient than the corresponding algorithms in ICMPv1 --- about 11 times faster for gA and at least 30 times faster for others. 

Finally, we  did tests on an ion channel molecular structure with two  ion channel pores to test the robustness of ICMPv2. We created this test case  through rotating  the $\alpha$-helix  of VDAC  by $20^{\circ}$ along the $z$ axis at the hinge region {\em (Gly-21–Tyr-22–Gly-23–Phe-24–Gly-25)}. A view of the molecular structure of  this modified VDAC, denoted by mVDAC,  is displayed in Figure~\ref{fig:3emnRotateZ20Comparison}.  With the values of parameters $Z1, Z2$, and $h_m$ and the box domain $\Omega$ used in the case of VDAC (see Table~1), ICMPv2 produced both membrane and solvent meshes in about 4.1 seconds only, showing the efficiency of our new schemes and the robustness of ICMPv2. A view of the protein, membrane, and solvent meshes of mVDAC is displayed in Figure~\ref{fig:3emnRotateZ20Comparison}. Here, the TMSmesh parameters $d=0.4$ and $c=e=0.9$; the numbers of vertices are 80630,16001,121943, 55779, 29185, and 66988 and the numbers of tetrahedra are 413086, 72646, 752618, 274675,138411, 339532, respectively, for the meshes $\hat{D}_{s,h}$,  $D_{s,h}^p$, $\Omega_h$, $D_{s,h}$, $D_{m,h}$, and  $D_{p,h}$.

As comparison, we did tests on mVDAC using ICMPv1 too. However, ICMPv1 was found not to work on this test case. The best meshes that we produced from ICMPv1 (in the sense of containing a small number of false tetrahedra)  were displayed in Figure~\ref{fig:3emnRotateZ20Comparison}. From Plot (c) it can be seen that the membrane mesh $D_{m,h}$ still contains many tetrahedra that belong to the solvent mesh $D_{s,h}$. Thus, the solvent mesh $D_{s,h}$ losses many tetrahedra so that its geometric shape has been twisted. Clearly, such poor membrane and solvent meshes may affect the approximation accuracy of a PBE/PNP finite element solution. Hence, this test case indicates that updating ICMPv1 to ICMPv2 is important and necessary.

\section{Conclusions}

A PBE/PNP ion channel finite element solver can effectively handle the complicated geometries of a protein region, a membrane region, and a solvent region as well as an interface between two of these three regions within a simulation box domain. This remarkable feature make it a much more powerful ion channel simulation tool than the corresponding finite difference solver. However, its numerical accuracy strongly depends on the  quality of an unstructured tetrahedral mesh to be used in its computer implementation. The application of a PBE/PNP ion channel finite element solver can be greatly enhanced and extended through developing and improving meshing algorithms and software. To further improve an ion channel mesh software package developed in Lu's research group a few years ago, which is referred as ICMPv1,  in this paper, based on ICMPv1, we have presented three new numerical schemes and implemented them. This work has resulted in the second version of ICMPv1, denoted as ICMPv2 here. Numerical tests on four ion channel proteins with different geometric complexities were done in this work, confirming that ICMPv2 can significantly improve the mesh quality and computer performance of ICMPv1 in the generation of an expanded solvent mesh and in the extraction of membrane and solvent meshes from a given expanded solvent mesh. Even so, more numerical tests are required to do to further explore the robustness and performance of  ICMPv2 for more complex ion channel proteins. We plan to do so in the future, along with the application of  unstructured tetrahedral meshes generated from ICMPv2 to the development of  PBE/PNP ion channel finite element solvers.

\section*{Acknowledgements}
This work was partially supported by  the Simons Foundation, USA,  through research award 711776. 
Gui and Lu's work was supported by the National Natural Science Foundation of China (Grant numbers 11771435 and 22073110).


\end{document}